\documentclass[hidelinks,onefignum,onetabnum]{siamart250211}

\usepackage{lipsum}
\usepackage{amsfonts}
\usepackage{graphicx}
\usepackage{epstopdf}
\ifpdf
  \DeclareGraphicsExtensions{.eps,.pdf,.png,.jpg}
\else
  \DeclareGraphicsExtensions{.eps}
\fi

\newsiamremark{remark}{Remark}
\newsiamremark{hypothesis}{Hypothesis}
\crefname{hypothesis}{Hypothesis}{Hypotheses}
\newsiamthm{claim}{Claim}
\newsiamremark{fact}{Fact}
\crefname{fact}{Fact}{Facts}

\headers{DualTPD for Nonlinear PDEs}{L. Chen, R. Guo, J. Wei, and J. Zou}

\title{A Novel Preconditioning Framework for Solving Nonlinear PDEs based on Fenchel-Rockafellar Duality and Transformed Primal-Dual Techniques \thanks{Submitted to the editors DATE.
\funding{The work of the first author was funded in parts by NSF DMS-2309777, and DMS-2309785. The work of the second author was funded by NSFC 12571436. The work of the third and fourth author was substantially supported by Hong Kong RGC General Research Fund (projects 14310324  and 14308322) and NSFC/Hong Kong RGC Joint Research Scheme (project N\_CUHK465/22). }}}

\author{
Long Chen \thanks{Department of Mathematics, University of California, Irvine, CA 92697, USA. 
 (\email{chenlong@math.uci.edu}).}
\and Ruchi Guo \thanks{School of Mathematics, Sichuan University, Chengdu, Sichuan, China, 610207
   (\email{ruchiguo@scu.edu.cn}).
  }
  \and Jingrong Wei \thanks{Department of Mathematics, The Chinese University of Hong Kong, Shatin, Hong Kong SAR, China
  (\email{jrwei@math.cuhk.edu.hk}).
  }
\and Jun Zou \thanks{Department of Mathematics, The Chinese University of Hong Kong, Shatin, Hong Kong SAR, China
  (\email{zou@math.cuhk.edu.hk}).}
}
\usepackage{amsopn}

\usepackage{amsfonts, psfrag, graphicx, mathtools, stackrel, amssymb,enumitem,indentfirst,float,color,enumerate,subcaption,pifont,mathtools,mathrsfs,caption}
\usepackage{epstopdf,wrapfig,chemarrow}

\usepackage{booktabs, cellspace, hhline}

\usepackage{algorithm}
\usepackage{multicol,multirow}
\usepackage{tcolorbox}
\usepackage{algpseudocode}

\usepackage{comment,multicol,multirow, xspace, diagbox,makecell}

\usepackage{tikz,tikz-cd}
\usepackage[utf8]{inputenc}
\usepackage{pgfplots} 
\usepackage{pgfgantt}
\usepackage{pdflscape}
\pgfplotsset{compat=newest} 
\pgfplotsset{plot coordinates/math parser=false}
\newlength\fwidth

\definecolor{myBlue}{rgb}{0.0,0.0,0.55}
\usepackage[hyperpageref]{backref}

\usepackage{comment,enumerate,xspace, diagbox}
\usepackage{booktabs, cellspace, hhline}

  \newcounter{mnote}
  \let\oldmarginpar\marginpar
    \renewcommand\marginpar[1]{\-\oldmarginpar[\raggedleft\footnotesize #1]%
    {\raggedright\footnotesize #1}}

\newtheorem{example}[theorem]{Example}
\newtheorem{algorithminline}[theorem]{Algorithm}

\newtheorem{assumption}[theorem]{Assumption}

\newcommand{\dd}{\,{\rm d}}

\newcommand{\curl}{{\rm curl\,}}
\renewcommand{\div}{\operatorname{div}}

\newcommand{\inprd}[1]{\left \langle#1 \right \rangle} % inner product

\newcommand{\mA}{\mathcal A}
\newcommand{\mB}{\mathcal B}

\newcommand{\mE}{\mathcal E}
\newcommand{\mF}{\mathcal F}
\newcommand{\mG}{\mathcal G}
\newcommand{\mI}{\mathcal I}

\newcommand{\mS}{\mathcal S}
\newcommand{\mT}{\mathcal T}
\newcommand{\mU}{\mathcal U}
\newcommand{\mV}{\mathcal V}

\newcommand{\bfC}{{\bf C}}

\newcommand{\bfD}{{\bf D}}

\newcommand{\bff}{{\bf f}}
\newcommand{\bfG}{{\bf G}}

\newcommand{\bfH}{{\bf H}}

\newcommand{\bfI}{{\bf I}}

\newcommand{\bfJ}{{\bf J}}

\newcommand{\bfM}{{\bf M}}

\newcommand{\bfR}{{\bf R}}
\newcommand{\bfr}{{\bf r}}
\newcommand{\bfS}{{\bf S}}

\newcommand{\bfu}{{\bf u}}

\newcommand{\bfphi}{\boldsymbol{\phi}}

\newcommand{\bfsigma}{\boldsymbol{\sigma}}

\ifpdf
\hypersetup{
  pdftitle={A Novel Preconditioning Framework for Solving Nonlinear PDEs based on Fenchel-Rockafellar Duality and Transformed Primal-Dual Techniques},
  pdfauthor={Long Chen, Ruchi Guo, Jingrong Wei, and Jun Zou}
}
\fi

\begin{document}

\maketitle

% REQUIRED
\begin{abstract}
A DualTPD method is proposed for solving nonlinear partial differential equations. The method is characterized by three main features. First, decoupling via Fenchel--Rockafellar duality is achieved, so that nonlinear terms are discretized by discontinuous finite element spaces, yielding block-diagonal mass matrices and closed-form updates. Second, improved convergence is obtained by applying transformed primal--dual (TPD) dynamics to the nonlinear saddle-point system, which yields strongly monotone behavior. Third, efficient preconditioners are designed for the elliptic-type Schur complement arising from the separated differential operators, and multigrid solvers are applied effectively. Extensive numerical experiments on elliptic $p$-Laplacian and nonlinear $H(\curl)$ problems are presented, showing significant efficiency gains with global, mesh-independent convergence.
\end{abstract}

% REQUIRED

% REQUIRED

\section{Introduction}

%Given an open bounded polyhedral Lipschitz domain $\Omega  \subset  \mathbb{R}^n$, we consider the convex minimization problem:
%\begin{equation}\label{eq: intro convex min}
%      \inf_{u \in \mV} \mI(u(x)) := \int_{\Omega} \mF(Du) + \mG(u) \dd x,
%\end{equation}
%where $D$ is a linear differential operator, $\mV$ is a Sobolev space such that $Du$ is well-defined in a certain sense, and $\mF, \mG$ are convex and differentiable energy densities. This formulation is widely used in partial differential equations (PDEs), for example, to represent the sum of kinetic and potential energy. The PDE formulation is derived by taking the first-order derivative of \eqref{eq: intro convex min}, leading to the Euler--Lagrange (E--L) equation. 

% In the optimization literature, \eqref{eq: min functional} is also called a composite problem (with a linear transformation), which appears widely in image recovery, signal processing, and partial differential equations (PDEs)~\cite{2007CombettesPesquet,goldstein2009split,nesterov2013gradient}. 

%In particular, the gradient (Frechét derivative) of $\mF$ at any fixed $\gamma \in \Sigma$ along a direction $\delta \in \Sigma$ can be represented as $\inprd{ \nabla \mF(\gamma), \delta }_{\Sigma}$, i.e., $ \nabla \mF(\gamma) \in \Sigma^*$.
%Henceforth, we shall use $\gamma = Du$ as an auxiliary variable. In various examples, $u$ and $\gamma$ may represent different variables; see Section \ref{sec: examples} for details. 

We consider the functional minimization problem
\begin{equation}\label{eq: min functional}
     \inf_{u \in \mV} \mI(u(x)) := \mE(u, Du) =  \mF(Du) + \mG(u), 
\end{equation}
where $\Sigma$ and $\mathcal{V}$ are reflexive Banach spaces, $D: \mV \rightarrow \Sigma$ is a linear operator, and $\mF: \Sigma \rightarrow \mathbb{R}$, $\mG: \mV \rightarrow \mathbb{R}$ are convex and differentiable. We develop fast optimization algorithms for solving problems of this form.

When $D$ is a differential operator, the formulation in \eqref{eq: min functional} can represent in physics the sum of kinetic and potential energy. Its minimization corresponds to solving the Euler--Lagrange (E--L) equation:
\begin{equation}\label{eq: E-L eq}
     \inprd{ \nabla \mF(D u), Dv} + \inprd{\nabla \mG(u), v} = 0, \quad \forall ~ v \in \mV,
\end{equation}
where $\inprd{\cdot, \cdot}$ denotes the duality pairing in the corresponding Banach space. The E--L equation \eqref{eq: E-L eq} leads to various nonlinear PDEs, such as the $p$-Laplacian problems~\cite{barrett1993finite, huang2007preconditioned} and nonlinear ferromagnetism~\cite{hichmani2025mixed} studied in this work. The well-posedness of \eqref{eq: min functional} and \eqref{eq: E-L eq} follows from the classical direct method in the calculus of variations; see Section \ref{sec: PD form} for details.

\begin{figure}[h]
\centering
\includegraphics[width=5.05in]{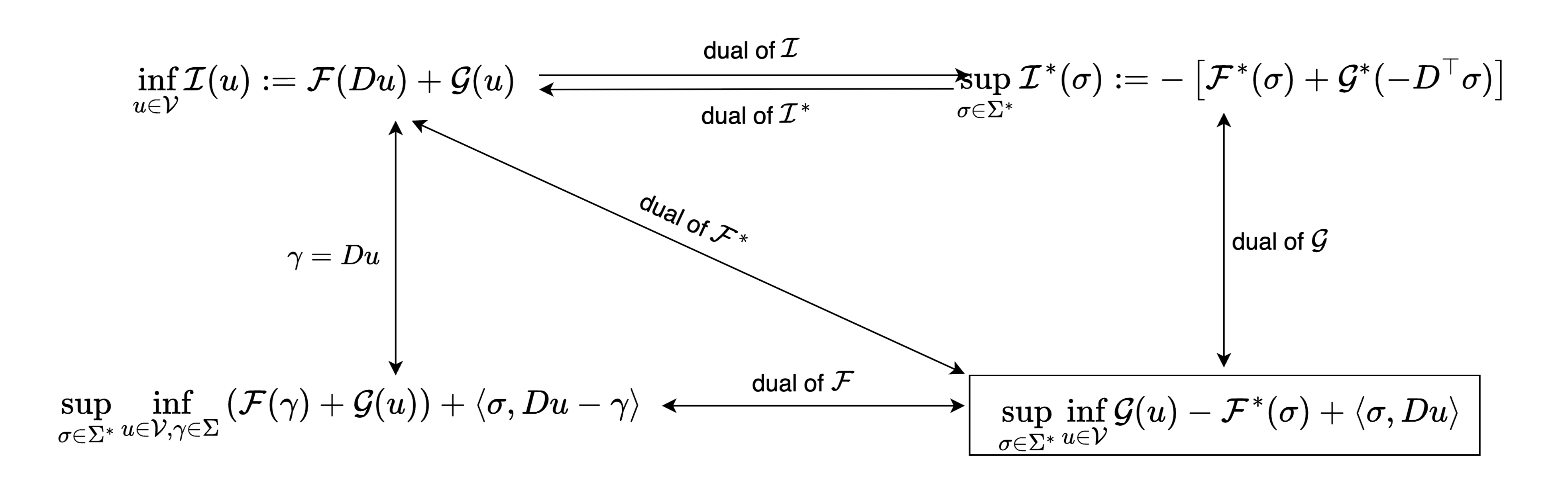}
\caption{The relation of various problems derived via Fenchel–Rockafellar duality.}
\label{fig_diagram}
\end{figure}

Developing fast optimization algorithms remains a major challenge in solving \eqref{eq: min functional} or \eqref{eq: E-L eq} directly. 
In numerical PDEs, notable approaches include the full approximation scheme~\cite{brandt1977multi,brune2015composing}, preconditioned descent algorithms~\cite{huang2007preconditioned}, and Newton-type methods~\cite{neuberger2013newton,pollock2015regularized}.
Fixed point methods~\cite{potter1993dual} are also widely used for \eqref{eq: E-L eq} because of their simple implementation.
In general, most of these methods target the Euler--Lagrange equation \eqref{eq: E-L eq},
that is, the upper-left block of the diagram in Figure \ref{fig_diagram},
where the main effort lies in designing efficient preconditioners. Using Fenchel–Rockafellar duality, we can connect \eqref{eq: min functional} to different formulations as shown in \cref{fig_diagram}. The focus of this work is to develop fast algorithms based on the sup--inf problem in the bottom-right block. 

%treating the nonlinearity $\mF$ and the differential operator $D$ simultaneously.

%To the authors' opinion, the difficulty primarily stems from two sources: the differential operator $D$ and the nonlinearity $\mF$. 
%Since $D$ is typically an unbounded operator, the spectrum of the discretized system may extend to infinity, 
%often leading to severe ill-conditioning—even in the linear case. 
%The situation is further complicated by the nonlinearity of $\mF$, which can distort the spectrum distribution and exacerbate the convergence behavior of nonlinear solvers. Collectively, these factors significantly hinder the development of efficient optimization methods.
%This has motivated significant interest in decoupling methods,
%which offer many favorable features such as exact solutions of nonlinear subproblems~\cite{aragon2023effective}.

The main difficulty in solving \eqref{eq: min functional} or \eqref{eq: E-L eq} comes from two sources: 
the differential operator $D$, 
which causes ill-conditioning of the resulting discrete systems, and the nonlinear $\mF$, 
which further increases the condition number and slows convergence of iterative solvers. 
This motivates decoupling methods, which can, for example, allow exact solutions of nonlinear subproblems~\cite{chambolle2011first}.
A common approach introduces $\gamma = \mathcal{C} Du$ and enforces this relation with a Lagrange multiplier,
where $\mathcal{C}$ is a linear operator depending on the coefficients; 
see, e.g.,~\cite{fortin2000augmented,gabay1983chapter,gabay1976dual}. This corresponds to solving the saddle-point problem in the bottom-left block of Figure \ref{fig_diagram}.
In numerical PDEs, this framework leads to mixed finite element methods (FEMs)~\cite{aragon2023effective,creuse2007posteriori,farhloul2000mixed}. 
Many iterative schemes have been developed by introducing auxiliary variables to split the objective functions, such as forward--backward splitting~\cite{combettes2005signal}, Douglas--Rachford splitting~\cite{2007CombettesPesquet}, Bregman splitting~\cite{goldstein2009split}, 
and others~\cite{liu2023fast}.

In this work, we propose a framework for developing fast iterative solvers for \eqref{eq: min functional} based on the {\it Fenchel--Rockafellar duality}~\cite{borwein2006convex,rockafellar1997convex,ekeland1999convex}.
Let $\mF^*$ (resp., $\mG^*$) denote the convex conjugate (Legendre transform) of $\mF$ (resp., $\mG$).
Suppose the Sobolev space $\Sigma$ is chosen so that the distributional dual pairing $\inprd{D^{\top}\sigma, u} := \inprd{\sigma, Du}$ is well defined.
Then the dual problem of \eqref{eq: min functional} is 
\begin{equation}\label{eq: intro dual_eq1}
     \sup_{\sigma\in \Sigma^*}  \mI^*(\sigma) := -  \left[ \mF^*(\sigma)+\mG^*(-D^{\top}\sigma)\right], 
\end{equation}
which is equivalent to the sup--inf problem:
\begin{equation}\label{eq: intro dual_eq2}
    \sup_{\sigma\in \Sigma^*}  \inf_{u\in\mV} \mS(\sigma,u) := \mG(u) -  \mF^*(\sigma)  + \inprd{\sigma,  Du }.
\end{equation}
The weak duality $\mI^*(\sigma) \leq \mI(u)$ holds for any $u \in \mV$, $\sigma \in \Sigma^*$. 
For brevity, further details are deferred to the next section.
Notice that \eqref{eq: intro dual_eq2} can be also rewritten as a constrained optimization problem,
solved by projected gradient methods \cite{1960Rosen,2025GuoZou} and primal-dual-type to be discussed below.

The Fenchel--Rockafellar duality is widely used in scientific computing.
For instance, it has served as a theoretical tool in~\cite{bartels2021error,carstensen2015nonconforming,carstensen2021unstabilized,tran2024discrete} to provide sharp error estimates of FEMs, although the computations in these works are still carried out on the primal problems.
It has also been observed that the decoupling property of this technique can be exploited to design fast solvers.
Notable approaches include the alternating direction method of multipliers (ADMM)~\cite{fortin2000augmented}, the Uzawa algorithm~\cite{1958ArrowHurwiczUzawa}, the Gauss--Newton method~\cite{carstensen2012mixed}, 
and primal--dual (PD) methods~\cite{attouch2016strongly,bartels2021nonconforming,briceno2011monotone+,chambolle2011first,zosso2017efficient} for various PDE-related optimization problems. Compared with the saddle-point problem in the bottom-left block of Figure \ref{fig_diagram}, one variable is eliminated but the knowledge of $\mF^*$ is required. This often leads to solving nonlinear subproblems~\cite{bartels2021nonconforming} or exploiting special structures of proximal updates~\cite{chambolle2011first,zosso2017efficient}. 

In this work, we show that the duality technique \eqref{eq: intro dual_eq2}, when combined with the transformed primal--dual (TPD) method~\cite{chen2023transformed,chen2024transformed,2025ZengZhanGuoWei} and suitable mixed finite element spaces, yields fast optimization algorithms, referred to here as the DualTPD method. 
TPD is a modified version of primal--dual methods, admitting better dynamical properties;
see detailed discussion in Section \ref{sec: algorithm}.
The proposed DualTPD algorithm \eqref{eq: EE TP-TPD} has three main features:
\begin{itemize}[leftmargin=20pt]

\item[(i)] {\it Decoupling via Fenchel--Rockafellar duality:} the duality framework separates the differential and nonlinear operators, 
allowing each to be treated with specialized techniques.
With mixed FEMs, the nonlinear term can be discretized by discontinuous FE spaces,
leading to a nonlinear block-diagonal mass matrix, which allows closed-form Jacobians and inverses.

\item[(ii)] {\it Improved convergence via transformed primal--dual dynamics:} 
applying the TPD method to the nonlinear saddle-point system produces a dynamical system with better monotonicity properties, which improves convergence.

\item[(iii)] {\it Efficient preconditioners:} 
the separated differential operators yield an ill-posed, elliptic-type Schur complement.
With suitably designed multigrid-based preconditioners, this term can be treated efficiently.

\end{itemize}
%\smallskip

The nonlinear saddle-point system is traditionally viewed as more difficult to solve. However, our DualTPD method can outperform solvers based on the primal formulation in the case of strong nonlinearity. Compared with the TPD methods developed in~\cite{chen2023transformed,chen2024transformed}, the present variant exploits duality to decouple the nonlinearity from the differential operators. Compared with the primal--dual method for the dual formulation~\cite{attouch2016strongly,bartels2021nonconforming,briceno2011monotone+,chambolle2011first,zosso2017efficient}, the transformed primal--dual method achieves faster convergence. Compared with Newton iteration (an outer iteration for nonlinearity combined with an inner linear saddle-point solver), our DualTPD couples the linearization and the solver in a single loop and retains fast global convergence.

We test our algorithm on two types of nonlinear equations: elliptic $p$-Laplacian problems and nonlinear $H(\curl)$ problems, and verify global and mesh-independent convergence of our algorithms. A full convergence analysis of the proposed algorithm will be explored elsewhere.

The paper is organized as follows. Section \ref{sec: PD form} derives the dual formulation of nonlinear PDEs and presents two examples: the $p$-Laplacian and the nonlinear ferromagnetism model, both of which fit into our framework. Section \ref{sec: algorithm} introduces the DualTPD algorithm based on the duality framework and transformed primal--dual techniques. Sections \ref{sec:dual p-Laplacian} and \ref{sec:dual Maxwell} demonstrate the efficiency of the proposed algorithm for the $p$-Laplacian and the nonlinear ferromagnetism model with designated preconditioners. Finally, Section \ref{sec:conclusion} provides concluding remarks.

%\input{Preliminaries}

% !TEX root =  DualTPD.tex

\section{Primal and Dual Formulations}\label{sec: PD form}
In this section, we recall the basic settings of the primal and dual formulations for optimization problems at the abstract level and present applications to nonlinear PDEs. We assume $\mF$ is proper, differentiable, strictly convex, and $\mG$ is a bounded linear functional.

\subsection{The Fenchel--Rockafellar Duality}
Introduce a new variable $\gamma = Du$ and reformulate the minimization problem \eqref{eq: min functional} as the constrained optimization problem
\begin{equation}\label{eq: constrained opt}
\inf_{u \in \mV, \gamma \in \Sigma, \,\gamma = Du} \mI(u,\gamma) := \mE(u, \gamma).
\end{equation}
With a Lagrange multiplier $\sigma$ to enforce the constraint $\gamma = Du$, 
\eqref{eq: constrained opt} is equivalent to the min--max problem
\begin{equation}\label{eq: saddle1}
   \sup_{\sigma \in \Sigma^*} \inf_{u \in \mV, \gamma \in \Sigma} 
   \Big( \mF(\gamma) + \mG(u) + \inprd{\sigma, Du - \gamma} \Big).
\end{equation}

The convex conjugate (Legendre transform) of $\mF$~\cite{rockafellar1997convex} is
\begin{equation}\label{eq: convex conjugate}
    \mF^*(\sigma) := \sup_{\gamma \in \Sigma} \big( \inprd{\sigma, \gamma} - \mF(\gamma) \big).
\end{equation}
Since $\mF$ is strictly convex, the maximum in \eqref{eq: convex conjugate} is attained when
\begin{equation}\label{eq: dual relation}
    \sigma = \nabla \mF(\gamma), 
    \qquad \text{equivalently,} \qquad
    \gamma = \nabla \mF^*(\sigma),
\end{equation}
see~\cite[Theorem 23.5]{rockafellar1997convex}. One can formally write $\nabla \mF^* = \nabla \mF^{-1}$ and $\nabla \mF = (\nabla \mF^*)^{-1}$.

With the conjugate function, \eqref{eq: saddle1} is equivalent to the inf-sup problem \eqref{eq: intro dual_eq2},
where the Lagrange multiplier coincides with the dual variable $\sigma$ in \eqref{eq: dual relation}.  
For a saddle point $(u, \sigma)$ of \eqref{eq: saddle1}, 
the first-order optimality (KKT) conditions are
\begin{subequations}\label{eq_saddle_intro}
    \begin{align}
 \inprd{\nabla \mF^*(\sigma), \tau} - \inprd{Du, \tau} &= 0,  \quad \forall~ \tau \in \Sigma^*, \\
 \inprd{\sigma, Dv} + \inprd{\nabla \mG(u), v} &= 0, \quad \forall~ v \in \mV,
    \end{align}
\end{subequations}
which, from a physical perspective, corresponds to a Hamiltonian system when \eqref{eq: min functional} represents a Lagrangian in mechanics~\cite{landau1960mechanics,arnol2013mathematical}.
%For a saddle point $(u, \sigma)$ of \eqref{eq: saddle}, the KKT condition results in the saddle point system in \eqref{eq_saddle_intro}.
%the following KKT conditions hold as the first-order optimality condition:
%\begin{subequations}\label{eq: KKT}
%    \begin{align}
%  \label{eq:KKT1}  &\inprd{\nabla \mF^*(\sigma), \delta} - \inprd{Du, \delta} = 0,  \quad \forall~ \delta \in \Sigma, \\
%\label{eq:KKT2}   & \inprd{\sigma, Dv}  + \inprd{\nabla \mG(u), v}= 0, \quad \forall~ v \in \mV,
%    \end{align}
%\end{subequations}
%which is also known as the Hamiltonian system if \eqref{eq: min functional} represents certain Lagrangian mechanics~\cite{landau1960mechanics,arnol2013mathematical}.
%We shall consider the case that $\mG$ is a linear mapping, i.e., $\mG(v) = \inprd{ g, v}$.

Moreover, using the identity $\inprd{D^{\top} u, \sigma} := \inprd{u, D\sigma}$, we have
\begin{equation}\label{rem_FR_dual_eq1}
\begin{split}
\inf_{u \in \mV} \big( \mG(u) + \inprd{\sigma, Du} \big) 
&= - \sup_{u \in \mV} \left( \inprd{-D^{\top}\sigma, u} - \mG(u) \right)  \\
&= -\mG^*(-D^{\top}\sigma),
\end{split}
\end{equation}
which yields the equivalence to \eqref{eq: intro dual_eq1} (top-right in \cref{fig_diagram}). 
%\mnote{Is $\mG$ linear or a general convex function? What is $\mG^*$ when $\mG$ is linear?}
%
This process reflects the Fenchel--Rockafellar duality~\cite{rockafellar1997convex}.  

In summary, the primal and dual formulations are equivalent, as illustrated in Figure \ref{fig_diagram} (top-left and bottom-right). We refer to \eqref{eq: min functional} as the primal problem, whose optimality condition is the Euler--Lagrange equation \eqref{eq: E-L eq}, and to \eqref{eq: intro dual_eq2} as the dual problem, whose optimality condition is the KKT system \eqref{eq_saddle_intro}. Our goal is to develop fast solvers for the nonlinear saddle-point system \eqref{eq_saddle_intro}.

We close this section by summarizing the results regarding the well-posedness.
Results for the primal problem can be found in~\cite{1997Evans,ekeland1999convex}. 
Here, we focus on the dual problem and recall the argument from~\cite{1977Scheurer}.  
%\mnote{ now we may need to assume $\mG$ is linear.}
Throughout the paper, we impose the following assumptions on $\mF^*$ and $D$:

\begin{assumption}\label{assump: F^*}
The function $\mF^*$ satisfies either global strong convexity with local Lipschitz continuity, or local strong convexity with global H\"older continuity in the Banach space $\Sigma^*$. That is, one of the following conditions holds:
 \begin{itemize}
 \item[(A1)] There exist $p > 2$ and constants $\mu, L > 0$ such that for all $\sigma_1, \sigma_2 \in \Sigma^*$,
 \begin{equation*}
\begin{aligned}
   & \inprd{\nabla \mF^*(\sigma_1) - \nabla \mF^*(\sigma_2), \sigma_1 - \sigma_2 } \geq \mu \|\sigma_1 - \sigma_2\|_{\Sigma^*}^p, \\
   & \|\nabla \mF^*(\sigma_1) - \nabla \mF^*(\sigma_2)\|_{\Sigma} \leq L (\|\sigma_1\|_{\Sigma^*} + \|\sigma_2\|_{\Sigma^*})^{p-2} \|\sigma_1 - \sigma_2\|_{\Sigma^*}.
\end{aligned}
\end{equation*}

 \item[(A1$^{\prime}$)] There exist $1 < p < 2$ and constants $\mu, L > 0$ such that for all $\sigma_1, \sigma_2 \in \Sigma^*$,
 \begin{equation*}
\begin{aligned}
   & (\|\sigma_1\|_{\Sigma^*} + \|\sigma_2\|_{\Sigma^*})^{2-p} 
   \inprd{\nabla \mF^*(\sigma_1) - \nabla \mF^*(\sigma_2), \sigma_1 - \sigma_2 } \geq \mu \|\sigma_1 - \sigma_2\|_{\Sigma^*}^2, \\
   & \|\nabla \mF^*(\sigma_1) - \nabla \mF^*(\sigma_2)\|_{\Sigma} \leq L \|\sigma_1 - \sigma_2\|_{\Sigma^*}^{p-1}.
\end{aligned}
\end{equation*}
 \end{itemize}
\end{assumption}

\begin{remark}
    There are other (weaker) conditions one may consider; see, e.g., \cite{ekeland1999convex,1997Evans}, but this is not the main interest in this work.
\end{remark}

\begin{assumption}\label{assump: D inf-sup}
The linear operator $D$ satisfies the inf--sup condition:
 \begin{equation}\tag{A2}
\inf_{v \in \mV} \sup_{\sigma \in \Sigma^* \backslash \{0\}} 
\frac{ \inprd{Dv,\sigma}_{\Sigma} }{\| v \|_{\mV} \, \| \sigma \|_{\Sigma^*}} \ge \beta,
\end{equation}
for some $\beta \in (0,\infty)$. 
\end{assumption}

\begin{theorem}[{\cite[Proposition 2.2]{1977Scheurer}}]
\label{thm_wellposed}
Let $\mF^*$ be proper, differentiable, and convex. Under Assumptions \ref{assump: F^*} and \ref{assump: D inf-sup}, 
the KKT system \eqref{eq_saddle_intro}, associated with the dual problem, admits a unique solution. 
\end{theorem}

\begin{remark}
% The global strong convexity and H\"older continuity are equivalent to each other by the argument in~\cite{1995DominiqueJeanPaul}.
   When $\mF$ instead of $\mF^*$ is accessible, we can verify strong convexity and H\"older continuity with respect to $\mF$ and then employ the duality argument~\cite{1995DominiqueJeanPaul} to obtain \cref{assump: F^*}.
\end{remark}

\subsection{Examples}\label{sec: examples}
We present two examples of nonlinear PDEs that fit into our formulation. 
In the cases considered here, $\mF$ is typically defined as the integral of a nonlinear function,  so each $\sigma$ acts as a linear functional on $\mV$ and can be identified with an element of $\mV^*$. We keep this relation throughout the discussion.

\begin{example}[$p$-Laplacian]\label{ex: p-lap}
The elliptic $p$-Laplacian problem~\cite{breit2015finite,farhloul2000mixed,creuse2007posteriori} is a typical degenerate nonlinear system and presents many challenges for finite element approximation and numerical solvers. For $\Omega \subseteq \mathbb{R}^d$, $d = 2,3$, $p > 1$ and $p^* = p/(p-1)$, let $\mV = W^{1,p}_0(\Omega)$, $\Sigma = [L^{p}(\Omega)]^d$, and $D = \nabla$. We denote the duality pairing $\inprd{g,v} := \int_{\Omega} g v \, \dd x$ and $|\cdot|$ as the $l^2$ norm of a vector. The nonlinear functional is given by
\begin{equation}\label{eq:p-laplace1}
 \mF(\gamma) = \int_{\Omega} |\gamma|^p \, \dd x, 
 \qquad \gamma = \nabla u,
\end{equation}
with the derivative defined as 
\begin{equation}\label{eq_pLap_F}
\langle \nabla \mF(\gamma), \delta \rangle = \int_{\Omega} |\gamma|^{p-2} \gamma \, \delta \, \dd x, 
\quad \forall \, \delta \in L^p(\Omega).
\end{equation}

The Euler--Lagrange equation is to find $u \in W^{1,p}_0(\Omega)$ such that
\begin{equation}\label{eq: ex p-laplace}
\inprd{|\nabla u|^{p-2} \nabla u, \nabla v} = \inprd{g, v}, 
\quad \forall~ v \in W^{1,p}_0(\Omega).
\end{equation}

From \eqref{eq_pLap_F}, the dual variable is $\sigma = |\gamma|^{p-2}\gamma$, and it is straightforward to verify that
\[
\gamma = \nabla \mF^*(\sigma) = |\sigma|^{p^*-2} \sigma,
\qquad
\mF^*(\sigma) = \int_{\Omega} |\sigma|^{p^*} \, \dd x.
\]

%In Darcy--Forchheimer problems modeled by $p$-Laplacians with $p = 1.5$~\cite{}, 
%the constitutive relation between $\sigma$ (pressure) and $\nabla u$ (velocity) is nonlinear.
The dual problem of \eqref{eq: ex p-laplace} is then
\begin{equation}\label{plap_dual_eq0}
\begin{aligned}
\inprd{ |\sigma|^{p^* - 2} \sigma, \tau} - \inprd{\nabla u, \tau} &= 0, 
~~~~~~\quad \forall~ \tau \in [L^{p^*}(\Omega)]^d, \\
\inprd{\sigma, \nabla v} &= \inprd{f, v}, 
\quad \forall~ v \in W_0^{1,p}(\Omega).
\end{aligned}
\end{equation}
$\mF$ and $\mF^*$ satisfy (A1) or (A1$^{\prime}$) naturally, and we refer interested readers to the detailed proof in \cite{Lindqvist2019}. (A2) also holds by the choice of spaces. 
% \JW{More details to  be added as appendix.}
% \LC{Continue with the well-posedness. Verify the assumptions.}
\end{example}

%\begin{example}
%    Consider the nonlinear curl-curl problem~\cite{2005BachingerLangerSchberl,2020XuYouseptZou}: find $u \in  H_0(\curl;\Omega) \cap \ker(\div;\Omega)$ such that
%\begin{equation*}
%\begin{aligned}
%\inprd{\nu(|\curl u|) \curl u, \curl v } &= \inprd{\mathcal{J}, v}, ~~~~ \forall ~ v \in H_0(\curl;\Omega) \cap \ker(\div;\Omega).
%\end{aligned}
%\end{equation*}
%The given source term $\mathcal J \in L^2(\Omega)$ satisfies
%$\div \mathcal J = 0 $. In this case, $D = \curl$, $\nabla \mF(v) = \nu(|v|)v$ and $\nabla \mG(u) = \mathcal J$. 
%\end{example}

\begin{example}[Nonlinear ferromagnetism]
\label{ex_ferro}
The second example is the ferromagnetism model~\cite{bachinger2005numerical,xu2020adaptive}.
Here, $D = \curl$, 
and $u$ denotes the magnetic potential with the magnetic field $\gamma := Du$. 
%This example is more involved than Example \ref{ex: p-lap}, since Gauss's law imposes the constraint $\div u = 0$, 
%which we enforce with a Lagrange multiplier.
The nonlinearity is given by a positive function $\nu(t): \mathbb{R}_{\geq 0} \to \mathbb{R}_{\geq 0}$ such that 
$\Phi(t) := \nu(t)t$ is strictly monotone.
Then the potential is defined as
\begin{equation}\label{ex_ferro_eq1}
 \mathcal{F}(\gamma) = \int_{\Omega} \int_0^{|\gamma|} \nu(t) t \, \dd t \, \dd x,
 \qquad
 \gamma = \curl u.
\end{equation}

This example is more complicated than Example \ref{ex: p-lap} as there is a constraint $\div u = 0$ from Gauss's Law.
Let $\bfH(\div0;\Omega)$ denote the divergence-free space.
We take $\mV = \bfH_0(\curl;\Omega)\cap \bfH(\div0;\Omega)$ and 
$\Sigma = [L^2(\Omega)]^3$ as Hilbert spaces. 
In the numerical experiments in Section \ref{sec:dual Maxwell}, we also consider Banach spaces, which usually correspond to stronger nonlinearities. 
Compared with existing work that requires the mapping $s \mapsto \nu(s)s$ to be strongly monotone and Lipschitz continuous~\cite{yousept2013optimal}, 
we relax the condition to $\mF$ satisfying (A1) or (A1'), which allows a broader class of examples, including the $p$-curl problems. (A2) holds due to the Poincar\'e inequaltiy $\|v\|\lesssim \|\curl v\|$ for $v\in \bfH_0(\curl;\Omega)\cap \bfH(\div0;\Omega)$.

The Euler--Lagrange equation is to find the magnetic potential $u \in \mV$ such that
\begin{equation}\label{eq:ex maxwell}
\inprd{\nu(|\curl u|)\,\curl u, \curl v } = \inprd{g, v}, 
\qquad \forall ~ v \in \bfH_0(\curl;\Omega).
\end{equation}
%In this case, $D = \curl$, $\mV = \bfH_0(\curl;\Omega)$ and 
%$\Sigma = \left[ L^2(\Omega) \right]^3$. 
%The coefficient function $\nu(t): \mathbb{R}_{\geq 0} \to \mathbb{R}_{\geq 0}$ is a positive nonlinear function such that $\Phi(t):=\nu(t)t$ is strictly monotone.
%Thus, the inverse function $\Phi^{-1}(s)$ of $\nu(t)t$ is well-defined, i.e., 
%$$
%s = \Phi(t) = \nu(t)t \Leftrightarrow t = \Phi^{-1}(s).
%$$
%The potentials are defined as follows:
%\begin{equation}
%\begin{split}
%\label{ex_ferro_eq3}
% \mathcal{F}(\gamma) = \int_{\Omega} \int_0^{|\gamma|} \nu(t)t \dd t \dd x, ~~~~~~~~~~ \mF^*(\sigma) = \int_{\Omega} \int_0^{|\sigma|}\varphi(t) \dd t \dd x.  
%\end{split}
%\end{equation}
%The dual formulation of nonlinear Maxwell equations \eqref{ex_ferro_eq2} is nontrivial and we shall investigate in Section \ref{sec:dual Maxwell}.

Introduce the dual variable $\sigma = \nu(|\gamma|)\gamma$,
and one can verify that 
\begin{equation}\label{ex_ferro_eq2}
\gamma = \nabla \mF^*(\sigma) = \frac{1}{\nu(\Phi^{-1}(|\sigma|))}\sigma,
\qquad
\mF^*(\sigma) = \int_{\Omega} \int_0^{|\sigma|} \Phi^{-1}(t) \, \dd t \, \dd x.
\end{equation}
The nonlinear relation between $\sigma$ and $\gamma$, from a physical perspective, represents the $B$--$H$ curve, 
whose specific form depends on the material. 

%To impose the constraint $\div u = 0$, we introduce the Lagrange multiplier $\phi \in H_0^1(\Omega)$.
%The dual formulation of \eqref{eq:ex maxwell} is then to find $\sigma \in [L^p(\Omega)]^3$, $\phi \in H_0^1(\Omega)$, and $u \in \bfH_0(\curl;\Omega)$ such that
%\begin{equation}\label{maxwell_dual_eq0}
%\begin{aligned}
%\inprd{\tfrac{1}{\nu(\Phi^{-1}(|\sigma|))}\sigma, \tau} - \inprd{\curl u, \tau} &= 0, 
%\quad \forall ~ \tau \in [L^p(\Omega)]^3, \\
%\inprd{u, \nabla \xi} &= 0, 
%\quad \forall ~ \xi \in H^1_0(\Omega), \\
%\inprd{\sigma, \curl v} + \inprd{\nabla \phi, v} &= \inprd{g, v}, 
%\quad \forall ~ v \in \bfH_0(\curl;\Omega),
%\end{aligned}
%\end{equation}
The dual formulation of \eqref{eq:ex maxwell} is then to find $\sigma \in \Sigma$ and $u \in \mV$ such that
\begin{equation}\label{maxwell_dual_eq0}
\begin{aligned}
\inprd{\tfrac{1}{\nu(\Phi^{-1}(|\sigma|))}\sigma, \tau} - \inprd{\curl u, \tau} &= 0, 
\quad \forall ~ \tau \in [L^p(\Omega)]^3, \\
%\inprd{u, \nabla \xi} &= 0, \quad \forall ~ \xi \in H^1_0(\Omega), \\
\inprd{\sigma, \curl v} &= \inprd{g, v}, 
\quad \forall ~ v \in \bfH_0(\curl;\Omega),
\end{aligned}
\end{equation}
where the given source term $g \in L^2(\Omega)$ satisfies $\div g = 0$.
%It is straightforward to check that the multiplier satisfies $\phi = 0$. 
% \LC{Continue with the well-posedness. Verify the assumptions.}
\end{example}

%%%%%%%%%%%%%%%%%%%%%%%%%%%%%%%%%%%%%%%%%%%%%%%%%

\subsection{Discretization with Block-diagonal Mass Matrices}\label{sec: Block-diagonal Discretization}
In this subsection, we present the discrete system of the dual formulation (in Hamiltonian mechanics) by finite element spaces and show that the nonlinearity leads to block-diagonal mass matrices. 
This structure allows the nonlinear terms to be inverted analytically, which significantly simplifies implicit or Newton-type solvers. 

The PDE problems are posed on an open domain $\Omega \subset \mathbb{R}^d$ 
and modeled in Sobolev spaces $\Sigma$ and $\mV$. 
Let $\mathcal T_h = \bigcup_i T_i$ denote a shape-regular triangulation of $\Omega$. 
Since no differential operator is applied to $\sigma$ in \eqref{eq_saddle_intro}, 
the space $\Sigma^*$ requires only minimal regularity and can be discretized by $L^2$-subspaces $\Sigma^*_h$, 
for example discontinuous polynomial spaces on the mesh:
\[
\Sigma^*_h = \{\sigma_h \in \Sigma : \; \sigma_h|_{T} \in [P_k(T)]^d , ~ \forall T \in \mathcal T_h \}.
\]
By construction, Assumption \ref{assump: F^*} holds naturally on $\Sigma^*_h \subset \Sigma^*$. 
As this space is fully broken, with no inter-element continuity, 
the term $\inprd{\nabla \mF^*(\sigma), \tau}$ discretizes to a block-diagonal matrix, 
so all nonlinear operations can be performed elementwise. 
Moreover, as $\Sigma^*_h=\Sigma_h$ in the discrete level, we shall simply use $\Sigma_h$ in our discussion below.

The discretization space $\mV_h$ for $\mV$ is chosen so that $Dv_h$ is well defined for $v_h \in \mV_h$. 
In particular, we take
\begin{equation}\label{Vh}
\mV_h := \{ v_h \in \mV : \; Dv_h \in [P_{k}(T)]^d \},
\end{equation}
which must satisfy the discrete inf--sup condition corresponding to Assumption \ref{assump: D inf-sup}:
\begin{itemize}
\item[(A$2_h$)] $\mV_h$ satisfies the inf--sup property
 \begin{equation*}
\inf_{v_h \in \mV_h} \sup_{\sigma_h \in \Sigma^*_h \backslash \{0\}} 
\frac{ \inprd{Dv_h, \sigma_h} }{\| v_h \|_{\mV_h} \, \| \sigma_h \|_{L^2(\Omega)}} \ge \bar{\beta},
\end{equation*}
for some $\bar{\beta} > 0$ independent of the mesh size $h$.
\end{itemize}

With these discretization spaces, the discretized problem is to
find $\sigma_h\in \Sigma_h$ and $u_h\in \mV_h$ satisfying
\begin{subequations}\label{eq_saddle_dis}
    \begin{align}
 \inprd{\nabla \mF^*(\sigma_h), \tau_h} - \inprd{Du_h, \tau_h} &= 0, ~ \quad\quad\quad \forall~ \tau_h \in \Sigma_h, \\
 \inprd{\sigma_h, Dv_h} &= \inprd{g, v_h} , \quad \forall~ v_h \in \mV_h,
    \end{align}
\end{subequations}
where $\mG(v):= -\inprd{ g, v}$. Under Assumptions (A1) (or (A1$^{\prime}$)) and (A$2_h$), 
the well-posedness of the discretized problem \eqref{eq_saddle_dis} is guaranteed. 
We remark that, in order for $\mV_h$ to satisfy (A$2_h$) in the case of nonlinear Maxwell equations in \cref{ex_ferro},
we shall employ Lagrange multiplier to impose $\div u = 0$ in a weak sense, 
which results in a larger saddle point system;
see details in Section \ref{sec:dual Maxwell}.

%\LC{Discrete problem.}

 %%%%%%%%%%%%%%%%%%%%%%%%%%%%%%%%%%%%%%%%%%%%%%%%

\section{The DualTPD Algorithm}\label{sec: algorithm}
In this section, we derive an algorithm for solving the saddle-point system \eqref{eq_saddle_intro} or \eqref{eq_saddle_dis} arising from the dual formulation of PDEs. We first review the transformed primal--dual algorithm from~\cite{chen2024transformed} and then proposes an improved variant so that it fits for efficiently solving the dual system.

\subsection{Transformation and Preconditioning}

For clarity, we rewrite the nonlinear saddle-point system \eqref{eq_saddle_intro} in matrix form: 
find $(\sigma, u) \in \Sigma \times \mV$ such that
\begin{equation*}
    \begin{pmatrix}
       \nabla \mF^* & -D \\
        D^{\top} &  O
    \end{pmatrix}
    \begin{pmatrix}
        \sigma \\
        u
    \end{pmatrix} 
    =
    \begin{pmatrix}
        0 \\
        g
    \end{pmatrix},
\end{equation*}
where $\nabla \mF^* \sigma := \nabla \mF^*(\sigma)$. 
 
Given a linear operator $\mI_{\sigma}^{-1}$, we consider the transformation
\begin{equation}\label{eq: tranformed ssp}
 \underbrace{\begin{pmatrix}
       I & O \\
       -D^{\top} \mI_{\sigma}^{-1} & I
    \end{pmatrix}}_{:= \mT}     
 \underbrace{\begin{pmatrix}
       \nabla \mF^* & -D \\
        D^{\top} &  O
    \end{pmatrix}}_{:= \mA}
 \begin{pmatrix}
        \sigma \\
        u
    \end{pmatrix} 
    = 
    \mT \begin{pmatrix}
        0 \\
        g
    \end{pmatrix}
    =
    \begin{pmatrix}
        0 \\
        g
    \end{pmatrix}.
\end{equation}
Since $\mT$ is invertible, the solution of the transformed system \eqref{eq: tranformed ssp} is equivalent to that of the saddle-point system \eqref{eq_saddle_intro}. 
%and
%\begin{equation*}
%    \mT \mA = \begin{pmatrix}
%           \nabla \mF^*  & -D \\
%            D^* (I - \mI_{\sigma}^{-1} \nabla \mF^* )   &  D^*\mI_{\sigma} D
%    \end{pmatrix}.
%\end{equation*}    

In the linear case, when $\nabla \mF^*(\sigma) = A$ for a linear invertible operator $A$, choosing $\mI_{\sigma}^{-1} = A^{-1}$ gives
\[
\mT \mA = 
\begin{pmatrix}
           A & -D \\
           O & S
\end{pmatrix},
\]
an upper triangular system where $S = D^{\top} A^{-1} D$ is the Schur complement. 
In practice, however, $A^{-1}$ is not computable for general nonlinear problems.

%Therefore, the transformation gives us two messages:
%\begin{itemize}
%    \item[(i)] $\mI_{\sigma}^{-1}$ shall be a preconditioenr for $\nabla \mF^*(\sigma)$.
%
%    \item[(ii)] The transformed system $\mT \mA$ shall be preconditioned by an upper triangular preconditioner in the form of $\mU^{-1}:= \begin{pmatrix}
%       \mI_{\sigma}& -D \\
%        O & \mI_S
%    \end{pmatrix}^{-1} =  \begin{pmatrix}
%       \mI_{\sigma}^{-1}& \mI_{\sigma}^{-1} D \mI_S^{-1}\\
%        O & \mI_S^{-1}
%    \end{pmatrix} $ where $\mI_S^{-1}$ is the preconditioner for $D^{\top} \mI_{\sigma}D$.
%\end{itemize}

With the transformation, we consider the preconditioned system:
$$
\begin{aligned}
    \mB^{-1} \mA  := \mU^{-1} \mT \mA 
      &  = \begin{pmatrix}
       \mI_{\sigma}& -D \\
        O & \mI_u
    \end{pmatrix}^{-1}
    \begin{pmatrix}
       I & O \\
       -D^{\top} \mI_{\sigma}^{-1} &  I
    \end{pmatrix}
    \begin{pmatrix}
       \nabla \mF^* & -D \\
        D^{\top} &  O
    \end{pmatrix} \\
%    & = \left ( \begin{pmatrix}
%       I & O \\
%       D^{\top} \mI_{\sigma}^{-1} &  I
%    \end{pmatrix} \begin{pmatrix}
%       \mI_{\sigma}& -D \\
%        O & \mI_S
%    \end{pmatrix} \right )^{-1}\begin{pmatrix}
%       \nabla \mF^* & -D \\
%        D^{\top} &  O
%    \end{pmatrix} \\
%    & =   \begin{pmatrix}
%       \mI_{\sigma} & -D\\
%       D^{\top} & \mI_S - D^{\top} \mI_{\sigma}^{-1}D
%    \end{pmatrix}^{-1} \begin{pmatrix}
%       \nabla \mF^* & -D \\
%        D^{\top} &  O
%    \end{pmatrix} .
\end{aligned}
$$
%\begin{remark}
%The preconditioner $\mB^{-1}$ is generally a dense operator. 
With the decomposition $\mB^{-1} = \mU^{-1}\mT$, the computation only requires the block preconditioners $\mI_{\sigma}^{-1}$ and $\mI_u^{-1}$.
%\end{remark}
The dual formulation is crucial for the design of effective preconditioners $\mI_{\sigma}^{-1}$ and $\mI_u^{-1}$.

\subsection{Transformed Primal--Dual Algorithm for Dual System}
\label{subsec:tri_PD}
% Now we are ready to derive the algorithm. Based on the explicit Euler discretizaiton of the flow \eqref{eq: TP-TPD flow}, 
Based on the preconditioning, we propose the DualTPD algorithm for solving the saddle-point problem \eqref{eq_saddle_intro} or \eqref{eq_saddle_dis}: 

\begin{algorithminline}[DualTPD algorithm]
Given initial value $(\sigma_0, u_0)$, preconditioners $\{\mI_{\sigma_k}^{-1}\}_{k \geq 0}, \{\mI_{u_k}^{-1}\}_{k \geq 0}$ and step size $\alpha_k > 0$, the sequence $\{(\sigma_k, u_k)\}_{k \geq 0}$ is generated by

\begin{equation}\label{eq: EE TP-TPD}
    \begin{pmatrix}
    \sigma_{k+1} \\
    u_{k+1}
\end{pmatrix} 
= 
\begin{pmatrix}
    \sigma_k \\
    u_k
\end{pmatrix} 
- \alpha_k \mB_k^{-1}
\begin{pmatrix}
   \nabla \mF^*(\sigma_k ) - Du_k \\
    D^{\top} \sigma_k - g
\end{pmatrix},
\end{equation}
where 
\[
\begin{aligned}
    \mB_k^{-1}:=  \mU_k^{-1} \mT_k
       = \begin{pmatrix}
       \mI_{\sigma_k}& -D \\
        O & \mI_{u_k}
    \end{pmatrix}^{-1}
    \begin{pmatrix}
       I & O \\
       -D^{\top} \mI_{\sigma_k}^{-1} &  I
    \end{pmatrix}.
\end{aligned}
\]
\end{algorithminline}

The iteration is most practical in residual--correction form:

\textit{Step 1:} Given $(\sigma_k, u_k)$, compute the residuals:
\begin{equation*}
    r_k^{\sigma} = \nabla \mF^*(\sigma_k) - Du_k, 
    \qquad 
    r_k^u = D^{\top} \sigma_k - g.
\end{equation*}

\textit{Step 2:} Apply the transformation and preconditioners $\mI_{\sigma_k}^{-1}$ and $\mI_{u_k}^{-1}$ to compute the correction:
\begin{equation*}
    \delta u_k =  \mI_{u_k}^{-1}  \big(r_k^u - D^{\top} \mI_{\sigma_k}^{-1} r_k^{\sigma}\big),
    \qquad
    \delta \sigma_k =  \mI_{\sigma_k}^{-1} \big(r_k^{\sigma} + D\delta u_k\big).
\end{equation*}

\textit{Step 3:} Update $(\sigma_{k+1}, u_{k+1})$ with step size $\alpha_k > 0$:
\begin{equation*}
    \sigma_{k+1} = \sigma_k - \alpha_k \delta \sigma_k,
    \qquad 
    u_{k+1} = u_k - \alpha_k \delta u_k.
\end{equation*}

Compared with the TPD algorithm in~\cite{chen2024transformed,chen2023transformed}, a key difference lies in the derivation of the saddle-point system  \eqref{eq_saddle_intro} or \eqref{eq_saddle_dis} from the dual formulation of the PDEs. We demonstrate that the duality-based framework is essential for effective preconditioning and enhanced computational efficiency. In the primal formulation \eqref{eq: E-L eq}, the nonlinear operator and differential operator are inherently coupled, complicating the design of preconditioners capable of addressing both challenges simultaneously. As shown in~\cite{chen2024transformed}, the absence of effective preconditioning leads to mesh-dependent iteration counts, particularly in highly nonlinear regimes.

On the other hand, the dual formulation \eqref{eq_saddle_intro} or \eqref{eq_saddle_dis} introduces a new potential computational consideration: the evaluation of $\nabla \mF^*(\sigma_k)$. While general dual descent algorithms requires solving a nonlinear subproblem~\cite{balseiro2020dual,yu2011dual}, our framework enables efficient evaluation of $\nabla \mF^*$ by leveraging the elementwise localization of nonlinearity within the discretized space $\Sigma_h$. Furthermore, the decoupling of the differential operator and the structure of $\Sigma_h$ allow for the design of Jacobian-based preconditioners $\mI_{\sigma_k}^{-1}$ and $\mI_{u_k}^{-1}$, and incorporating Newton-type solvers. By combining these techniques, our proposed algorithm the proposed algorithm achieves global and mesh-independent convergence for a broad class of nonlinear PDEs. % \LC{please be consistent with the notation. the correction is denoted by $\delta$ in algorithms.}

\subsection{Connection to Other Algorithms}
Despite being derived from the dual formulation, our algorithm applies to general saddle-point systems. 
For linear saddle-point problems, triangular preconditioners were first introduced by Bramble and Pasciak~\cite{bramble1988preconditioning}
and have since been extensively developed; see~\cite{benzi2005numerical,benzi2011modified,chen2018convergence, zulehner2002analysis} and the references therein. 
The motivation comes from the row reduction of saddle-point systems, later formalized as a transformation and extended to nonlinear saddle-point problems~\cite{chen2023transformed}. 
A key distinction between linear and nonlinear saddle-point systems lies in the choice of the relaxation parameter (step size): 
for linear problems, it is determined by spectral analysis, while for nonlinear problems, it is guided by Lyapunov analysis. In \cite{chen2024transformed}, we showed through Lyapunov analysis that TPD with variable preconditioners yields a globally convergent scheme.

For nonlinear saddle-point problems, diagonal preconditioners are often considered within the framework of primal--dual methods~\cite{chambolle2011first,zosso2017efficient} 
or inexact Uzawa algorithms~\cite{chen1998preconditioned,song2019inexact,chen1998global,hu2001iterative, hu2002two,hu2006nonlinear}. 
% In~\cite{chen1998preconditioned}, Chen studied the SOR--Newton method and multistep inexact Uzawa method. 
A nonlinear preconditioner $\mI_{\sigma}$ for $\nabla \mF^*$ can be used if the resolvent or its proximal operator has a closed form~\cite{chambolle2011first,song2019inexact,zosso2017efficient}, 
or if the corresponding subproblem can be solved up to a prescribed accuracy~\cite{chen1998global,hu2002two,hu2006nonlinear}. The primal--dual hybrid gradient method (PDHG)~\cite{chambolle2011first} is effective when the proximal operators of the nonlinear terms are available, and it has been widely applied to imaging problems~\cite{valkonen2014primal,sidky2012convex,fang2014single}. However, without effective preconditioners, PDHG often requires thousands of iterations to converge, especially for PDE problems with poor conditioning~\cite{carrillo2022primal,clason2017primal}.  
Notice that $(\nabla \mF^*)^{-1} = \nabla \mF$. 
If we use $\sigma_k = \nabla \mF(Du_k)$ and substitute this into the equation for $u$, the scheme reduces to gradient descent for the primal problem, which is typically slow without an appropriate preconditioner. 
In this setting, $\mI_u$ should should serve as an effective preconditioner for the primal problem. Choosing a linear 
$\mI_{\sigma}$ corresponds to linearizing the nonlinearity, which often makes construction of  $\mI_u$ substantially easier.

Augmented Lagrangian methods (ALM)~\cite{hestenes1969multiplier,rockafellar1976augmented} are also widely used for constrained optimization. 
By adding an augmented term $\nabla \mF^* + \beta D D^{\top}$ for $\beta > 0$, the resulting system corresponds to the dual problem \eqref{eq: intro dual_eq1}. One can use a scaled identity as the Schur complement preconditioner when $\beta \gg 1$~\cite{chen2023transformed}. 
In this case, large $\beta$ alleviates the difficulty of designing a good $\mI_u$. 
However, constructing an efficient preconditioner for $\nabla \mF^* + \beta D D^{\top}$ with large $\beta$ remains challenging.  

In particular, if $\mI_{\sigma} = \nabla^2 \mF^*(\sigma_k)$, the Hessian of $\mF^*$ at $\sigma_k$, and $\mI_{u} = D^{\top} \mI_{\sigma}^{-1} D$, with step size $\alpha_k = 1$, method \eqref{eq: EE TP-TPD} reduces to Newton's method. 
It is well known that Newton's method converges locally: quadratic convergence is achieved only when the initial guess is sufficiently close to the true solution, and there is no guarantee of global convergence. While Newton's method is mainly considered for quasilinear or semilinear cases~\cite{neuberger2013newton,pollock2015regularized}, where the nonlinear term is essentially decoupled from the differential operator, we show that the dual formulation \eqref{eq_saddle_intro} and \eqref{eq_saddle_dis} decouple the nonlinear term from the differential operator. 
This allows Newton-type preconditioners to be adapted either in explicit closed form or by solving simple algebraic equations after finite element discretization. 
Most importantly, the computations are local: the associated matrices are elementwise independent, and the inverses involve only $d \times d$ blocks in the domain $\Omega \subset \mathbb{R}^d$.

\section{Application I: $p$-Laplacian Equations}\label{sec:dual p-Laplacian}

In this section, we apply the triangularly preconditioned DualTPD algorithm to the $p$-Laplacian problem. 
We refer readers to \cref{ex: p-lap} for background. 
We show that the design of preconditioners benefits from the dual formulation, 
which, combined with duality and the block-diagonal structure, significantly accelerates convergence.

\subsection{Discretization and Matrix Form}
For computation, we discretize the system using conforming FE spaces and use elementwise constant elements for $\sigma \in \Sigma_h = V^T_h(\Omega)$ and linear Lagrange elements with the corresponding boundary condition for $u \in \mathcal{V}_h = V^n_h(\Omega)$:
\begin{equation*}
     \begin{aligned}
         V^n_h(\Omega) & = \{ u \in C^0(\Omega) : u|_{\partial \Omega} = 0, \; u|_{T} \in \mathbb{P}_1(T), \; \forall T \in  \mathcal{T}_h \}, \\
         V^T_h(\Omega) & = \{ \sigma \in [L^2(\Omega)]^2 : \sigma|_{T} \in [\mathbb{P}_0(T)]^2, \; \forall~ T \in  \mathcal{T}_h \}.
     \end{aligned}
\end{equation*} 
The dimensions of $V^T_h(\Omega)$ and $V^n_h(\Omega)$ are denoted by $N_T$ and $N_n$, respectively. 
By choosing bases of $\Sigma_h$ and $\mV_h$, 
functions in these spaces can be identified with vectors.
We denote the following matrix representations of the related operators:

\begin{itemize}[leftmargin=16pt]
\item With a fixed weight $\nu$, define $\bfM^{\nu}_{\dd} \in \mathbb{R}^{N_{\dd}\times N_{\dd}}$ as the mass matrix on $V^n_h(\Omega)$ or $V^T_h(\Omega)$, for $\dd = n, T$. 
% In particular, if $\nu = 1$, $\bfM_{\dd}$ denotes the standard mass matrix. 

\smallskip

\item The matrix associated with the weak divergence operator is $\bfD$. 

\smallskip

\item For a fixed $\bfsigma \in \mathbb{R}^2$, define the coefficient $\gamma(\bfsigma) := |\bfsigma|^{p^* -2}$.
\end{itemize}
\smallskip
 
The corresponding matrix form of \eqref{eq_saddle_dis} is to find vectors $\bfsigma$ and $\bfu$ such that
\begin{equation*}
     \begin{pmatrix}
          \bfM_{T}^{\gamma(\bfsigma)} &  -\bfD \\
          \bfD^{\top} & \bf0
     \end{pmatrix}
     \begin{pmatrix}
          \bfsigma \\
          \bfu  
     \end{pmatrix} 
     = 
     \begin{pmatrix}
          \mathbf{0} \\
          \bff
     \end{pmatrix},
\end{equation*}
where for $T \in \mathcal T_h$ and $\bfsigma|_T = \bfsigma_T \in \mathbb{R}^2$, the nonlinear term has the block-diagonal structure
\begin{equation*}
    \bfM_{T}^{\gamma(\bfsigma)}= {\rm diag}\left \{     \begin{pmatrix}
        \gamma(\bfsigma_{T_i}) |T_i| & 0 \\
        0 & \gamma(\bfsigma_{T_i}) |T_i|
    \end{pmatrix}
\right \}_{i = 1}^{N_T}, 
%    \quad 
%    \bfM^{\gamma}_{T} = 
%    \begin{pmatrix}
%        \gamma(\bfsigma_T) |T| & 0 \\
%        0 & \gamma(\bfsigma_T) |T|
%    \end{pmatrix},
    \quad  \bfM_{T}^{\gamma(\bfsigma)} \bfsigma = \big(\gamma(\bfsigma_{T_i}) |T_i| \bfsigma_{T_i}\big)_{i = 1}^{N_T}.
\end{equation*}
%and 
%\[
%\]
%%%%%%%%%%%%%%%%%%%%%%%%%%%%%%%%%%%%%%%%%%%%%%%

\subsection{Preconditioning}

We now discuss the preconditioners $\bfI_{\sigma_k}^{-1}$ and $\bfI_{u_k}^{-1}$ in the algorithm \eqref{eq: EE TP-TPD}. 
For nonlinear problems, a natural choice for $\bfI_{\sigma_k}^{-1}$ is the mass matrix with nonlinear coefficients. 
In particular, $\bfM_{T}^{\gamma(\bfsigma)}$ is a diagonal matrix for any given $\bfsigma$, and its inverse can be computed directly. 
Note that $\bfM_{T}^{\gamma(\bfsigma)}$ or its inverse becomes degenerate when $\bfsigma|_{T} =\bf0$ for some $T \in \mathcal T_h $.  Following \cite{huang2007preconditioned}, we adopt the regularization 
\begin{equation}\label{eq: mass precon}
\gamma_{\lambda}(\bfsigma) = \begin{cases}
\gamma(\bfsigma) + \lambda & \quad  p^* >  2, \\
 (\gamma |\bfsigma| + \lambda)^{p^*-2}, 
    & \quad p^* < 2 \ \text{and}\  |\bfsigma| \leq \epsilon_0, \\
 \gamma(\bfsigma), & \quad \text{otherwise}
\end{cases}
\end{equation}
% \begin{equation}\label{eq: mass precon}
%     \bfI_{\sigma_k}^{-1} = \begin{cases}
%     \left( \bfM_{T}^{\gamma(\bfsigma_k) + \epsilon}\right)^{-1}, & \quad  p^* >  2, \\[6pt]
%     \left( \bfM_{T}^{\gamma_{\epsilon}(\bfsigma_k)}\right)^{-1}, 
%     \quad \gamma_{\epsilon}(\bfsigma_k)  = (|\bfsigma_k| + \epsilon)^{p^*-2}, 
%     & \quad p^* < 2 \ \text{and}\  |\bfsigma_k| \leq \epsilon_0, \\[6pt]
%      \left( \bfM_{T}^{\gamma(\bfsigma_k)}\right)^{-1}, & \quad \text{otherwise},
% \end{cases}
% \end{equation}
where $\lambda > 0$ is a regularization parameter to avoid singular matrix inversion, 
and $\epsilon_0 > 0$ is the threshold for applying regularization. With the regularized coefficient, we have the first candidate of preconditioner $\bfI_{\sigma_k}^{-1} = \left (\bfM_{T}^{\gamma_{\lambda}(\bfsigma_k)} \right)^{-1}$. This scheme is referred to as  {\it DualTPD-M}. 

Next, we consider a Newton-type preconditioner based on the Jacobian of the nonlinear term. We shall show the Jacobian of $\bfM_{T}^{\gamma(\bfsigma)}\bfsigma$ preserves the block structure.
\begin{equation}\label{eq: pLap Jacobian}
    \bfJ(\bfM_{T}^{\gamma(\bfsigma)}\bfsigma) = {\rm diag}\{\bfJ_{T_i}\}_i, 
\end{equation}
where $\bfJ_T \in \mathbb{R}^{2\times 2}$ is given by
\begin{equation*}
\begin{aligned}
    \bfJ_T &= \gamma(\bfsigma_T) |T| \bfI_2 + |T| \bfsigma_T \nabla \gamma(\bfsigma_T)^{\top} \\
           &= |\bfsigma_T|^{p^*-2}|T| \bfI_2  + (p^* - 2)|\bfsigma_T|^{p^* - 4}|T|\bfsigma_T\bfsigma_T^{\top}.
\end{aligned}
\end{equation*}
% $\bfJ_T$ or its inverse is degenerate if and only if $|\bfsigma_T| = 0$. 
When $|\bfsigma_T| \neq 0$, this special form yields an obvious advantage that its inverse has a closed form by using the Sherman–Morrison formula:
\begin{equation*}
    \bfJ_T^{-1} = \frac{|\bfsigma_T|^{2 - p^*}}{|T|} \bfI_2 - \frac{(p^* - 2)|\bfsigma_T|^{-p^*}}{(p^* - 1)|T|} \bfsigma_T \bfsigma_T^{\top},
\end{equation*}
which shows the inverse of Jacobian \eqref{eq: pLap Jacobian}:
\[
 [\bfJ(\bfM_{T}^{\gamma(\bfsigma)}\bfsigma)]^{-1} =  {\rm diag}\{\bfJ_{T_i}^{-1}\}_i.
\]

As with $\gamma^{-1}(\bfsigma_T)$ and $|\bfsigma_T|^{-1}$, the inverse $\bfJ_T^{-1}$ is only applied when both are greater than the threshold $\epsilon_0$. 
More precisely, our second preconditioner is defined as follows: for a given $\bfsigma$, $\bfI_{\sigma}^{-1} = {\rm diag}\{\bfI_{T_i}^{-1}\}_{i=1}^{N_T}$ where $\bfI_{T}^{-1} \in \mathbb{R}^{2\times 2}$ is defined as
\begin{equation}\label{eq: Jacobian precon}
\bfI_{T}^{-1} = \begin{cases}
\begin{pmatrix}
        \gamma_{\lambda}(\bfsigma_{T}) |T| & 0 \\
        0 & \gamma(\bfsigma_{T}) |T|
    \end{pmatrix}^{-1}, & \quad  \text{ if } p^* > 2, \gamma(\bfsigma_{T}) \leq \epsilon_0 \text{ or } p^* < 2, |\bfsigma_T| \leq \epsilon_0,  \\[6pt]
  \bfJ_{T}^{-1}, & \quad \text{otherwise}.
\end{cases}
\end{equation}
This scheme is referred to as \textit{DualTPD-J}.  
% \begin{equation}\label{eq: Jacobian precon}
% \bfI_{\sigma_k}^{-1} = \begin{cases}
%     \left( \bfM_{T}^{\gamma(\bfsigma_k) + \epsilon}\right)^{-1}, & \quad  p^* > 2 \ \text{and}\ \gamma(\bfsigma_k) \leq \epsilon_0,  \\[6pt]
%     \left( \bfM_{T}^{\gamma_{\epsilon}(\bfsigma_k)}\right)^{-1}, 
%     \quad \gamma_{\epsilon}(\bfsigma_k)  = (|\bfsigma_k| + \epsilon)^{p^*-2}, 
%     & \quad p^* < 2 \ \text{and}\  |\bfsigma_k| \leq \epsilon_0, \\[6pt]
%     [\bfJ(\bfM_{T}^{\gamma(\bfsigma_k)}\bfsigma_k)]^{-1}, & \quad \text{otherwise}.
% \end{cases}
% \end{equation}

With such $\bfI_{\sigma_k}^{-1}$, the corresponding Schur complement is
\[
\bfS_k = \bfD^{\top} \bfI_{\sigma_k}^{-1}\bfD,
\]
which is a variable Laplacian matrix with tensor coefficients. 
We employ a geometric multigrid V-cycle as the preconditioner, i.e., 
\[
\bfI_{\sigma_k}^{-1} = {\rm MG}(\bfS_k, {\rm tol}_{\rm mg}, {\rm maxIt}),
\]
where the number of V-cycles is bounded by \rm{maxIt} or the relative residual is reduced below ${\rm tol}_{\rm mg}$.

If $\alpha_k = 1$ and $\bfI_{\sigma_k}^{-1} = \bfS_k^{-1}$, the method reduces to Newton's iteration.  

\subsection{Numerical Results}

In all experiments, we set ${\rm tol}_{\rm mg} = 10^{-2}$ and maxIt $= 5$ for the inner iteration. For the preconditioner, we set the regularization parameters $\lambda = 10^{-4}$ and $\epsilon_0 = 10^{-16}$. 
The algorithm terminates once the relative residual falls below $10^{-6}$:
\begin{equation*}
    {\rm rel}_{r} = \frac{|(\bfr^{\sigma}_k, \bfr^u_k)|}{|\bff|} \leq 10^{-6},
\end{equation*}
where 
\[
(\bfr^{\sigma}_k, \bfr^u_k) = \big( \bfM_{T}^{\gamma(\bfsigma_k)}\bfsigma_k - \bfD \bfu_k, \; \bfD^{\top} \bfsigma_k - \bff \big).
\] 
%For cases where the nonlinear coefficient is used to construct preconditioners as in \eqref{eq: mass precon}, 
%we refer to the algorithm as {\it DualTPD-M}. 
%For cases where the Jacobian is used to construct preconditioners as in \eqref{eq: Jacobian precon}, 
%we refer to the algorithm as {\it DualTPD-J}. 
We set the regularization parameters $\lambda = 10^{-4}$ and $\epsilon_0 = 10^{-16}$.

\smallskip

\begin{table}[ht]
  \centering
  \caption{ Numerical errors on $p$-Laplacian with $p = 4$.}
  \label{table:p-Laplacian-numerical-error}
  \renewcommand{\arraystretch}{1.25}
    \resizebox{6.5cm}{!}{  
  \begin{tabular}{@{} c c c c @{}}
    \toprule
$h$ & DoF	&	$\|u - u_h\|_{L^2}$	&	$\|\sigma - \sigma_h\|_{L^2}$	\\ \hline
1/32 & 12417	&	1.45E-05	&	2.11E-04	\\
1/64 & 49409	&	3.61E-06	&	1.06E-04	\\
1/128 & 197121	&	9.04E-07	&	5.29E-05	\\
1/256 & 787457	&	2.26E-07	&	2.64E-05	\\
1/512 & 3147777	&	5.65E-08	&	1.32E-05	\\
    \bottomrule
  \end{tabular}}
\end{table}

For the first test, we validate the convergence of the finite element approximation. 
Let $\Omega = [0,1]^2$ with the exact solution $u(x,y) = 10x(x-1)y(y-1)$. 
Table~\ref{table:p-Laplacian-numerical-error} reports the numerical errors with respect to the mesh size $h$ and degrees of freedom (DoF) for $p=4$. The observed convergence rates are $\mathcal{O}(h^2)$ for the primal variable $u$ and $\mathcal{O}(h)$ for the dual variable $\sigma$, verifying the correctness of our finite element discretization. 
We next examine the performance of the solvers. 

We consider both $p < 2$ and $p > 2$ cases with zero initial values. 
Table~\ref{table:p-Laplacian-numerical-result-ex1} reports the iteration counts, averaged multigrid V-cycles per iteration, and CPU time for our proposed algorithms with $p = 1.5$ and $p = 4$. 
Step sizes were first determined through preliminary runs on coarse meshes and then fixed for all mesh sizes $h$. 

%For $p = 1.5$, the step size is $\alpha_k = 1$ for DualTPD-J and $\alpha_k = 0.6$ for DualTPD-M.  
%For $p = 4$, the step size is $\alpha_k = 0.6$ for DualTPD-J and $\alpha_k = 1.2$ for DualTPD-M.  

% \begin{table}[htp]
%   \centering
%   \caption{ Computation results on $p$-Laplacian with $p = 1.5$ and $p =4$.}
%   \label{table:p-Laplacian-numerical-result-ex1}
%   \renewcommand{\arraystretch}{1.15}
%     \resizebox{11.75cm}{!}{  
%   \begin{tabular}{@{} c | c || c c || c c @{}}
%     \toprule
%     &		&	\multicolumn{2}{c||}{$p = 1.5$}	&	\multicolumn{2}{c}{$p = 4$}	\\ \cline{3-6}
% 	&		&	DualTPD-J	&	DualTPD-M	&	DualTPD-J	&	DualTPD-M	\\ 
% 	&	$h \backslash\alpha$	&	1	&	0.6	&	0.6	&	1.2	\\ \hline
% \multirowcell{ 5}{Iteration \\ (MG V-cycle)}	&	1/32	&	5 (25)	&	16 (80)	&	20 (90)	&	26 (104)	\\
% 	&	1/64	&	5 (25)	&	17 (85)	&	20 (93)	&	 25 (115)	\\
% 	&	1/128	&	5 (25)	&	17 (85)	&	20 (95)	&	25 (118)	\\
% 	&	1/256	&	5 (25)	&	17 (85)	&	20 (98)	&	25 (118)	\\
% 	&	1/512	&	5 (25)	&	17 (85)	&	22 (110)	&	25 (118)	\\ \hline
% \multirowcell{ 5}{CPU time \\ in seconds}	&	1/32	&	0.043	&	0.11	&	0.18	&	0.17	\\
% 	&	1/64	&	0.17	&	0.38	&	0.6	&	0.62	\\
% 	&	1/128	&	0.46	&	1.1	&	1.8	&	1.6	\\
% 	&	1/256	&	1.6	&	3.9	&	7.1	&	5.6	\\
% 	&	1/512	&	7.3	&	16	&	34	&	24	\\
%     \bottomrule
%   \end{tabular}
%   }
% \end{table}

\begin{table}[htp]
  \centering
  \caption{ Computation results for the proposed DualTPD algorithm on $p$-Laplacian with $p = 1.5$ and $p =4$.}
  \label{table:p-Laplacian-numerical-result-ex1}
  \renewcommand{\arraystretch}{1.15}
    \resizebox{12.05cm}{!}{  
  \begin{tabular}{@{} c | c || c c || c c @{}}
    \toprule
    &		&	\multicolumn{2}{c||}{$p = 1.5$}	&	\multicolumn{2}{c}{$p = 4$}	\\ \cline{3-6}
	&		&	DualTPD-J	&	DualTPD-M	&	DualTPD-J	&	DualTPD-M	\\ 
	&	$h \backslash\alpha$	&	1	&	0.6	&	0.6	&	1.2	\\ \hline
\multirowcell{ 5}{Iteration \\ (Avg. MG per step)}	&	1/32	&	5 (5)	&	16 (5)	&	20 (4.5)	&	26 (4)	\\
	&	1/64	&	5 (5)	&	17 (5)	&	20 (4.7)	&	 25 (4.6)	\\
	&	1/128	&	5 (5)	&	17 (5)	&	20 (4.8)	&	25 (4.7)	\\
	&	1/256	&	5 (5)	&	17 (5)	&	20 (4.9)	&	25 (4.7)	\\
	&	1/512	&	5 (5)	&	17 (5)	&	22 (5)	&	25 (4.7)	\\ \hline
\multirowcell{ 5}{CPU time \\ in seconds}	&	1/32	&	0.043	&	0.11	&	0.18	&	0.17	\\
	&	1/64	&	0.17	&	0.38	&	0.6	&	0.62	\\
	&	1/128	&	0.46	&	1.1	&	1.8	&	1.6	\\
	&	1/256	&	1.6	&	3.9	&	7.1	&	5.6	\\
	&	1/512	&	7.3	&	16	&	34	&	24	\\
    \bottomrule
  \end{tabular}
  }
\end{table}

The results show that the iteration count is independent of the mesh size $h$. The CPU time grows nearly linearly with the number of degrees of freedom. 
Since each residual calculation involves a sparse matrix--vector product of size $\mathcal{O}(N)$, our proposed algorithm \eqref{eq: EE TP-TPD} achieves optimal computational scaling. For $p = 1.5$, as the zero initial value is close to the solution, DualTPD-J with $\alpha = 1$ converges in just $5$ iterations, exhibiting superlinear convergence and reflecting the local behavior of Newton’s method.

%\begin{figure}
%    \centering
%    \includegraphics[width=0.5\linewidth]{Figures/timegrowth_plap.png}
%    \caption{Rate of time growth for DualTPD algorithms to solve $p$-Laplacian.}
%    \label{fig:timegrowth_plap}
%\end{figure}

%
%\begin{figure}
%    \centering
%    \includegraphics[width=0.45\linewidth]{Figures/u_pllap.png}
%    \includegraphics[width=0.45\linewidth]{Figures/sigma_pllap.png}
%    \caption{Order of convergence for the numerical error by the presented algorithm on p-Laplacian.}
%    \label{fig:conver order p-Laplacian}
%\end{figure}
%\vspace{0.2cm}

Our second test is the nonhomogeneous Dirichlet problem on the unit disk $B(0,1)$ with $f \equiv 1$. 
The exact solution is 
\[
u(x) = \frac{p-1}{p}\left(\tfrac{1}{2}\right)^{\tfrac{1}{p-1}}\left(1-|x|^{\tfrac{p}{p-1}}\right), 
\quad  x \in B(0,1).
\]

We compare our proposed algorithm \eqref{eq: EE TP-TPD} with the primal problem solved by the preconditioned gradient descent (PGD) algorithm~\cite{huang2007preconditioned}. 
The parameters are set as suggested in that work. 
Since the most time-consuming parts of PGD are the line search and exact multigrid solves, we also implement a variant with fixed step size and inexact multigrid solves (denoted PGD w/ LS). 

\begin{table}[ht]
  \centering
  \caption{Iteration counts and CPU time for various algorithms on the $p$-Laplacian with $p = 1.5$ using zero and random initial values.}
  \label{table:p-Laplacian-comparison-merged}
  \renewcommand{\arraystretch}{1.15}
  \resizebox{11.75cm}{!}{    
  \begin{tabular}{@{} c | c | cc | cc | cc | cc @{}}
    \toprule
    & & \multicolumn{2}{c|}{DualTPD-J} & \multicolumn{2}{c|}{DualTPD-M} & \multicolumn{2}{c|}{PGD} & \multicolumn{2}{c}{PGD w/ LS} \\ \cline{3-10}
$p = 1.5$    & $h \backslash \alpha$ & \multicolumn{2}{c|}{1} & \multicolumn{2}{c|}{0.8} & \multicolumn{2}{c|}{line search} & \multicolumn{2}{c}{0.2} \\ %\cline{3-10}
    &  & Zero & Rand & Zero & Rand & Zero & Rand & Zero & Rand \\ \hline
    \multirow{3}{*}{Iteration} 
      & 1/32  & 4  & 11 & 12 & 23 & 43 & 51  & 134 & 172 \\
      & 1/64  & 4  & 12 & 11 & 23 & 38 & 37  & 133 & 182 \\
      & 1/128 & 4  & 14 & 10 & 24 & 56 & 60  & 133 & 192 \\ \hline
    \multirowcell{3}{CPU time \\ in seconds} 
      & 1/32  & 0.09 & 0.23 & 0.20 & 0.37 & 2.8 & 3.6 & 2.0 & 3.0 \\
      & 1/64  & 0.25 & 0.77 & 0.57 & 1.0  & 7.5 & 7.3 & 6.6 & 9.1 \\
      & 1/128 & 0.85 & 3.3  & 1.58 & 3.7  & 44  & 48  & 23  & 40 \\ 
    \bottomrule
  \end{tabular}
  }
\end{table}

\begin{table}[ht]
  \centering
  \caption{Iteration counts and CPU time for various algorithms on the $p$-Laplacian with $p = 10$ using zero and random initial values.}
  \label{table:p-Laplacian-comparison-2-merged}
  \renewcommand{\arraystretch}{1.15}
  \resizebox{11.75cm}{!}{    
  \begin{tabular}{@{} c | c | cc | cc | cc | cc @{}}
    \toprule
    & & \multicolumn{2}{c|}{DualTPD-J} & \multicolumn{2}{c|}{DualTPD-M} & \multicolumn{2}{c|}{PGD} & \multicolumn{2}{c}{PGD w/ LS} \\ \cline{3-10}
 $p = 10$    & $h \backslash \alpha$ & \multicolumn{2}{c|}{0.2} & \multicolumn{2}{c|}{1.2} & \multicolumn{2}{c|}{line search} & \multicolumn{2}{c}{0.2} \\ %\cline{3-10}
    &  & Zero & Rand & Zero & Rand & Zero & Rand & Zero & Rand \\ \hline
    \multirow{3}{*}{Iteration} 
      & 1/32  & 105 & 114 & 65 & 59 & 38 & 70  & 93 & 122 \\
      & 1/64  & 114 & 127 & 65 & 60 & 35 & 65  & 93 & 114 \\
      & 1/128 & 125 & 136 & 64 & 62 & 38 & 69  & 94 & 112 \\ \hline
    \multirowcell{3}{CPU time \\ in seconds} 
      & 1/32  & 2.2 & 2.4 & 1.1 & 0.97 & 2.5 & 4.1 & 1.5 & 2.0 \\
      & 1/64  & 7.7 & 8.3 & 2.9 & 2.8  & 6.2 & 12  & 4.5 & 5.6 \\
      & 1/128 & 28  & 31  & 9.3 & 9.1  & 27  & 47  & 16  & 19  \\ 
    \bottomrule
  \end{tabular}
  }
\end{table}

Tables~\ref{table:p-Laplacian-comparison-merged} and \ref{table:p-Laplacian-comparison-2-merged} report the iteration counts and CPU time for different methods with $p = 1.5$ and $p = 10$, using both zero and random initial values. 
In all cases, our proposed algorithms are more time-efficient compared with PGD since no line search is required and inexact multigrid is used as a preconditioner. For the PGD, these relaxations lead to a significant increase in the number of iterations. The results also demonstrate that our algorithms are robust to the choice of initial values and converge globally. 

% For each $p$, thechange of initial guess would 
% In particular with random initial guesses, our proposed algorithms remain more efficient than the primal algorithms for both $p < 2$ and $p > 2$ cases. 

%%%%%
In particular, our proposed algorithms handle the cases $p < 2$ efficiently, which are known to be unstable for the primal problem as $p \to 1$. 
To further support this claim, we report in \cref{table:p-Laplacian-Dual-merged} the iteration counts of our algorithms for a range of $p$ values. 
The regularization parameters are set to $\epsilon_0 = \lambda = 10^{-4}$ for $p < 2$ and $\lambda = 10^{-4}$, $\epsilon_0 = 10^{-16}$ for $p > 2$. 

We observe that the iteration counts and algorithm parameters are independent of the mesh size $h$. 
Moreover, DualTPD-M performs better for $p > 2$, while the Jacobian preconditioner DualTPD-J is more effective for $p < 2$. 

\begin{table}[htp]
  \centering
  \caption{Iteration numbers of DualTPD-J (Jacobian preconditioner) and DualTPD-M (mass preconditioner) for the $p$-Laplacian with corresponding step sizes $\alpha$.}
  \label{table:p-Laplacian-Dual-merged}
  \renewcommand{\arraystretch}{1.25}
  \resizebox{7.8cm}{!}{  
  \begin{tabular}{@{} c | c || c c c c c @{}}
    \toprule
    &   $p$  & 1.05 & 1.3 & 1.5 & 4 & 10 \\ 
    \hline
    \multirowcell{5}{DualTPD-J \\ Iteration \\} 
        & $h \backslash \alpha$ & 1   & 1   & 1   & 0.6 & 0.2 \\ \cline{2-7}
        & 1/32   & 5 & 6  & 4  & 22  & 105 \\
        & 1/64   & 5 & 6  & 4  & 22  & 114 \\
        & 1/128  & 5 & 6  & 4  & 22  & 125 \\
        & 1/256  & 5 & 6  & 4  & 22  & 127 \\ \hline
    \multirowcell{5}{DualTPD-M \\ Iteration \\} 
        & $h \backslash \alpha$ & 1   & 0.5 & 0.8 & 1.3 & 1.5 \\ \cline{2-7}
        & 1/32   & 5 & 22 & 12 & 17  & 65 \\
        & 1/64   & 5 & 22 & 11 & 17  & 65 \\
        & 1/128  & 5 & 22 & 10 & 17  & 64 \\
        & 1/256  & 5 & 23 & 10 & 17  & 64 \\
    \bottomrule
  \end{tabular}
  }
\end{table}

%If one were to always set $\alpha = 1$ for DualTPD-J, the iteration counts would likely change, but global convergence would no longer be guaranteed. 

%\begin{table}[ht]
%  \centering
%  \caption{Iteration number for the proposed algorithm on $p$-Laplacian with precondtioner \eqref{eq: mass precon}.}
%  \label{table:p-Laplacian-DualMass-p}
%  \renewcommand{\arraystretch}{1.25}
%  \begin{tabular}{@{} c | c c c c c @{}}
%    \toprule
%$h$	&	p = 1.05	&	p = 1.3	&	p = 1.5	&	p = 4	&	p = 10	\\ \hline
%1/32	&	5	&	22	&	12	&	17	&	65	\\
%1/64	&	5	&	22	&	11	&	17	&	65	\\
%1/128	&	5	&	22	&	10	&	17	&	64	\\
%1/256	&	5	&	23	&	10	&	17	&	64	\\
%    \bottomrule
%  \end{tabular}
%  \vspace{0.2cm}
%\end{table}
%
%
%%\begin{minipage}[t]{0.95\linewidth}
%\begin{table}[ht]
%  \centering
%  \caption{Iteration number for the proposed algorithm on $p$-Laplacian with precondtioner \eqref{eq: Jacobian precon}.}
%  \label{table:p-Laplacian-DualJacobian-p}
%  \renewcommand{\arraystretch}{1.25}
%  \begin{tabular}{@{} c | c c c c c @{}}
%    \toprule
% $h$	&	p = 1.05	&	p = 1.3	&	p = 1.5	&	p = 4	&	p = 10	\\ \hline
%1/32	&	5	&	6	&	4	&	22	&	105	\\
%1/64	&	5	&	6	&	4	&	22	&	114	\\
%1/128	&	5	&	6	&	4	&	22	&	125	\\
%1/256	&	5	&	6	&	4	&	22	&	127	\\
%    \bottomrule
%  \end{tabular}
%  \vspace{0.2cm}
%\end{table}

%  \end{minipage}

% !TEX root =  DualTPD.tex

\section{Application II: Nonlinear Maxwell Equations}\label{sec:dual Maxwell}

In this section, we study the dual formulation of nonlinear Maxwell equations. 
With suitable preconditioners, we demonstrate that the DualTPD algorithm \eqref{eq: EE TP-TPD} can be effectively applied to solve these problems.

\subsection{Augmented Lagrangian formulation}

Based on \eqref{maxwell_dual_eq0}, 
to enforce the divergence-free condition $\div u = 0$, 
we introduce a Lagrange multiplier $\phi$. 
To facilitate the construction of suitable preconditioners, 
we adopt the augmented Lagrangian technique by adding the term $\inprd{\phi,\xi}$, noting that the solution remains unchanged since $\phi = 0$. 
Thus, we consider the following modified dual problem, obtained from \eqref{maxwell_dual_eq0}: find $\sigma \in [L^2(\Omega)]^3$, $\phi \in H_0^1(\Omega)$, and $u \in \bfH_0(\curl;\Omega)$ such that
\begin{subequations}\label{eq:dual maxwell2}
\begin{align}
\inprd{\tfrac{1}{\nu(\Phi^{-1}(|\sigma|))}\sigma, \tau} - \inprd{\curl u, \tau} &= 0, \quad \forall ~ \tau \in [L^2(\Omega)]^3, \\
\label{eq:dual maxwell2_p} 
\beta \inprd{\phi,\xi} - \inprd{u, \nabla \xi} &= 0, \quad \forall ~ \xi \in H^1_0(\Omega), \\
\inprd{\sigma, \curl v} + \inprd{\nabla \phi, v} &= \inprd{\mathcal{J}, v}, \quad \forall ~ v \in \bfH_0(\curl;\Omega),
\end{align}
\end{subequations}
where $\beta > 0$ is a parameter to be specified later. 

Although \eqref{eq:dual maxwell2} involves three variables, 
it still fits into the general saddle-point structure \eqref{eq_saddle_intro} if $(\sigma, \phi)$ is treated as a single variable, 
with $\mF^*$ and $D$ modified accordingly. 
Since $\sigma$ and $\phi$ are decoupled, and $\phi$ enters only through linear terms, 
the extension of the algorithm \eqref{eq: EE TP-TPD} and the discussion of preconditioners carry over directly. 
% \breakline

% \begin{theorem}
%      Suppose \cref{assump: monotone nu} is satisfied. The solution to equation \eqref{eq:dual maxwell2} exists and is unique.
% \end{theorem}

% \begin{proof}
% \cref{lem: nabla f property} showed that $\nabla \mF \in \mS_{\underline{\nu}, L}$. As a result, the operator associated to $\nabla \mF^*(\sigma) := (\nabla \mF)^{-1}(\sigma)$ satisfies $\nabla \mF^* \in \mS_{L^{-1}, \underline{\nu}^{-1}}$ with $L^{-1} > 0$.

% Next, we introduce the weak divergence operator as the adjoint of negative gradient operator, i.e.,
% $$\inprd{\div_w u, p} : = - \inprd{ u, \nabla p},  \quad  u \in H_0(\curl; \Omega), p \in H_0^1(\Omega).$$

% According to dual formulation \eqref{eq:dual maxwell2}, $D = [-\curl, -\div_w]^{\top}$ and  
% \begin{equation*}
%     \begin{aligned}
%        \inprd{D u, Dv }   = \inprd{\curl u, \curl v} + \inprd{ \div_w u , \div_w v}, \quad  u, v \in \mV.
%     \end{aligned}
% \end{equation*}

% Since $D$ is injective \cite{xxx}, we verify the assumptions in \cref{thm: general well-posedness g}.
    
% \end{proof}

% \breakline

\subsection{Discretization and Matrix Form}

For computation, we use element-wise constant elements for $\sigma$, linear Lagrange elements for $\phi$, and lowest-order N\'ed\'elec edge elements for $u$: 
\begin{equation*}
     \begin{aligned}
          V^T_h(\Omega) &= \{\sigma \in [L^2(\Omega)]^3 : \; \sigma|_{T} \in [\mathbb{P}_0(T)]^3, \; \forall T \in  \mathcal{T}_h \}, \\
         V^n_h(\Omega) &= \{\phi \in C^0(\Omega) : \; \phi|_{\partial \Omega} = 0, \; \phi|_{T} \in \mathbb{P}_1(T), \; \forall T \in  \mathcal{T}_h \}, \\
         V^e_h(\Omega) &= \{ u \in \bfH_0(\curl; \Omega) : \; u|_{K} = a + b \times x, \; a, b \in \mathbb{R}^3 , \; \forall T \in \mathcal{T}_h\}.
     \end{aligned}
\end{equation*}

With these spaces, the discretized weak form of \eqref{eq:dual maxwell2} is:  
find $(\sigma_h, \phi_h, u_h) \in V^T_h(\Omega) \times V^n_h(\Omega) \times V^e_h(\Omega)$ such that
\begin{equation}\label{maxwell_dual_eq1}
\begin{aligned}
\inprd{\tfrac{1}{\nu(\Phi^{-1}(|\sigma_h|))}\sigma_h, \tau_h} - \inprd{\curl u_h, \tau_h} &= 0, 
\quad \forall \tau_h \in V^T_h(\Omega), \\
\beta \inprd{\phi_h,\xi_h} + \inprd{u_h, \nabla \xi_h} &= 0, 
\quad \forall \xi_h \in V^n_h(\Omega), \\
\inprd{\sigma_h, \curl v_h} + \inprd{\nabla \phi_h, v_h} &= \inprd{g, v_h}, 
\quad \forall v_h \in V^e_h(\Omega).
\end{aligned}
\end{equation}

By selecting $\mV_h  = V^e_h(\Omega) \cap \{ \inprd{v_h,\nabla p_h} = 0, ~ \forall p_h \in V^n_h(\Omega) \}$ and $\Sigma_h = V^T_h(\Omega)$, one can verify that \eqref{eq_saddle_dis} and \eqref{maxwell_dual_eq1} are equivalent. 
However, formulation \eqref{maxwell_dual_eq1} is more favorable in computation, as it directly enforces the divergence-free condition and facilitates the design of preconditioners (see the next subsection).

Choosing bases for these spaces identifies all functions with vectors. 
Denote the dimensions of $V^T_h(\Omega)$, $V^n_h(\Omega)$, and $V^e_h(\Omega)$ by $N_T$, $N_n$, and $N_e$, respectively. 
The associated matrix operators are defined as follows:  

\begin{itemize}[leftmargin=16pt]

\item For a weight function $\nu$, let $\bfM^{\nu}_{\dd}\in \mathbb{R}^{N_{\dd}\times N_{\dd}}$ denote the mass matrix on $V^n_h(\Omega)$, $V^e_h(\Omega)$, or $V^T_h(\Omega)$ for $\dd = n, e, T$. 
In particular, if $\nu=1$, $\bfM_{\dd}$ is the standard mass matrix.  

\smallskip

\item Let $\bfG:\mathbb{R}^{N_n} \rightarrow \mathbb{R}^{N_T}$ be the grad matrix, 
let $\bfC:\mathbb{R}^{N_e} \rightarrow \mathbb{R}^{N_T}$ be the curl matrix, 
%Those matrices are indeed the incidence matrices between sub-simplexes if the correct scaled bases are chosen. 
and let $\bfD = \bfG^\top \bfM_e$ be the matrix associated with the weak divergence operator. 
Denote the matrix $\tilde \bfC = \bfM_{T}\bfC $.

%\smallskip
%\item  The matrix associated with the week divergence operator is $\bfD := \bfG^{\top} \bfM_e$. The matrix $\tilde \bfC = \bfM_{T}\bfC $.

%\item Let $\bfC : \mathbb{R}^{N_e} \to \mathbb{R}^{N_t}$ be the curl matrix, and let $\bfD$ denote the matrix associated with the weak divergence operator.  
%Define $\tilde \bfC = \bfM_t \bfC$.  

\smallskip

\item For a fixed $\bfsigma \in \mathbb{R}^3$, define $\gamma(\bfsigma) := 1/\nu(\Phi^{-1}(|\bfsigma|))$.  

\end{itemize}
% The nonlinear PDE associated to the primal formulation reads as:
% \begin{equation}\label{eq: primal matrix}
%    ( \bfC^{\top}\bfM_{3T}^{\nu(|\bfC \bfu |)}\bfC +\bfB^{\top} \bfM^{-1}_e \bfB) \bfu  = \bfJ.
% \end{equation}
%Then, the matrix form for the nonlinear saddle point system is 
%=======

With these definitions, the nonlinear saddle-point system takes the matrix form
\begin{equation}
    \label{eq: saddle pt mat}
   \begin{pmatrix}
          \bfM_{T}^{\gamma(\bfsigma)} & \bf0 & -\tilde \bfC \\
          \bf0   &  \bfM_n^{\beta} & -\bfD \\
          \tilde \bfC^{\top} &  \bfD^{\top} & \bf0
    \end{pmatrix}
    \begin{pmatrix}
       \bfsigma \\
       \bfphi \\
       \bfu  
    \end{pmatrix} = 
\begin{pmatrix}
   \mathbf{0} \\
   \mathbf{0} \\
   \bfJ
\end{pmatrix},
\end{equation}
where for $T \in \mathcal T_h$ and $\bfsigma|_T = \bfsigma_T \in \mathbb{R}^3$, the nonlinear term has the block-diagonal structure
\begin{equation*}
    \bfM_{T}^{\gamma(\bfsigma)}= {\rm diag}\left \{     \begin{pmatrix}
        \gamma(\bfsigma_{T_i}) |T_i| &  &\\
         & \gamma(\bfsigma_{T_i}) |T_i| &\\
         & &  \gamma(\bfsigma_{T_i}) |T_i| 
     \end{pmatrix}
\right \}_{i = 1}^{N_t}, 
%    \quad 
%    \bfM^{\gamma}_{T} = 
%    \begin{pmatrix}
%        \gamma(\bfsigma_T) |T| & 0 \\
%        0 & \gamma(\bfsigma_T) |T|
%    \end{pmatrix},
    \quad  
\end{equation*}
and 
\begin{equation*}
    \bfM_{T}^{\gamma(\bfsigma)} \bfsigma = \big(\gamma(\bfsigma_{T_i}) |T_i| \bfsigma_{T_i}\big)_{i = 1}^{N_t}.
\end{equation*}

\subsection{Preconditioning}\label{sec: preconditioning}

%Based on the discretization of the dual variable $\bf\sigma$, we construct preconditioners using the Jacobian. 
In this case, we need to construct the preconditioners for the three variables: $\bfsigma$, $\bfphi$ and $\bfu$.
First, for $\bfsigma$, since it is piecewise constant, the Jacobian $\bfJ(\bfM_{T}^{\gamma(\bfsigma)}\bfsigma)$ is block diagonal with $3\times 3$ blocks on each element. We denote $\bfsigma_T = \bfsigma|_{T} \in \mathbb{R}^3$.  Direct computation using $\gamma(\bfsigma)  = 1/(\nu(\Phi^{-1}(|\bfsigma|))$ gives that
\begin{equation*}
    \bfJ(\bfM_{T}^{\gamma(\bfsigma)}\bfsigma) = {\rm diag}\{\bfJ_{T_i}\}_i, 
\end{equation*}
where $\bfJ_T \in \mathbb{R}^{3\times 3}$ is given by
\begin{equation*}
\begin{aligned}
    \bfJ_T &= \gamma(\bfsigma_T) |T| \bfI_3 + |T| \bfsigma_T \nabla \gamma(\bfsigma_T)^{\top} = \gamma(\bfsigma_T) |T| \bfI_3  - t(\bfsigma)|T|\bfsigma_T\bfsigma_T^{\top},
\end{aligned}
\end{equation*}
where 
$$
t(\bfsigma) =  \frac{\nu^{\prime}(\Phi^{-1}(|\bfsigma|)) }{\nu^{2}(\Phi^{-1}(|\bfsigma|)|\bfsigma|} (\Phi^{-1})^{\prime}(|\bfsigma|).
$$

% $$\bfJ(\bfM_{3T}^{\gamma(\bfsigma)}\bfsigma) =\bfM_{3T}^{\gamma(\bfsigma)-t\bfsigma\bfsigma^{\top}} $$
% is a weighted mass matrix with tensor coefficient $\gamma(\bfsigma)-t\bfsigma\bfsigma^{\top}$ on each element and
% $$t =  \frac{\nu^{\prime}(\Phi^{-1}(|\bfsigma|)) }{\nu^{2}(\Phi^{-1}(|\bfsigma|)|\bfsigma|} (\Phi^{-1})^{\prime}(|\bfsigma|).$$
With the inverse function theorem and $\Phi(x) = \nu(x)x$, we have 
$$(\Phi^{-1})^{\prime}(|\bfsigma|) = \frac{1}{\Phi^{\prime}(z)} = \frac{1}{\nu^{\prime}(z)z +\nu(z)},$$
where $z =\Phi^{-1}(|\bfsigma|)$. Notice $\Phi(\cdot)$ is a smooth function, we can apply Newton's method on scalar problems to solve for $z$ and the coefficient
$$
t(\bfsigma, z) =  \frac{\nu^{\prime}(z) }{\Phi^{\prime}(z)\nu^{2}(z)|\bfsigma|} .
$$

Using Woodbury matrix identity, we can compute
$$
\bfJ_T^{-1} =  \frac{\nu(z)}{|T|}\bfI_3  - \frac{t(\bfsigma, z) \nu^2(z)}{ (t(\bfsigma, z) \nu(z)|\bfsigma|^2 - 1) |T|} \bfsigma \bfsigma^{\top}  \in \mathbb{R}^{3 \times 3}
$$
and the preconditioner for $\bfsigma$
$$\bfI_{\sigma}^{-1}  = \big(\bfJ(\bfM_{T}^{\gamma(\bfsigma)}\bfsigma)\big)^{-1} = {\rm diag}\{\bfJ_{T_i}^{-1}\}_{i = 1}^{N_t}.$$

% $ \bfI_{\sigma}^{-1}  = (\bfJ(\bfM_{3T}^{\gamma(\bfsigma)}\bfsigma))^{-1}$, which is the inverse of the mass matrix $\bfM_{3T}$ equipped with tensor coefficient \\ $\left( \nu(z)\bfI - \frac{\nu(z)}{-1/(t\nu(z)) + |\bfsigma|^2} \bfsigma \bfsigma^{\top} \right)$.

% Following the discussion in Section \ref{sec: Block-diagonal Discretization}, we set 
% \[
% \bfI_{\sigma}^{-1} = \big(\bfJ(\bfM_{T}^{\gamma(\bfsigma)}\bfsigma)\big)^{-1},
% \]
% which is the inverse of the mass matrix $\bfM_t$ with the corresponding tensor coefficients.

For the multiplier variable $\bfphi$, we can apply the mass lumping so that $\bfM_n$ becomes diagonal, and simply take $\bfI_{\phi}^{-1} = (\bfM^{\beta}_n)^{-1}$.

Next, for $\bfu$, we consider the Schur complement
\[
\bfS(\bfsigma) = \tilde \bfC^{\top}\big(\bfJ(\bfM_{T}^{\gamma(\bfsigma)}\bfsigma)\big)^{-1}\tilde \bfC 
+ \bfD^{\top}(\bfM_n^{\beta})^{-1}\bfD.
\]
For any given $\bfsigma$, the first term is a $\curl$--$\curl$ matrix with block-diagonal tensor coefficients, and the second term is a $\div$--$\div$ matrix. 
We use the preconditioned conjugate gradient (PCG) method and denote the preconditioner by
\[
\bfI_u^{-1} := \text{PCG}(\bfS(\bfsigma), {\rm tol}_{\rm pcg}),
\]
where ${\rm tol}_{\rm pcg}$ is the tolerance for terminating the PCG iterations.
For the inner preconditioner within the PCG, we adopt a multigrid method with a modified Hiptmair--Xu preconditioner~\cite{hiptmair2007nodal}. 
Specifically, we use the Hiptmair--Xu decomposition with variational auxiliary matrices~\cite{kolev2009parallel} to approximate $(\bfS(\bfsigma))^{-1}$:
\[
(\bfS(\bfsigma))^{-1} \approx \bfR + \bfG \mB_{\rm grad}\bfG^{\top} + \boldsymbol{\Pi} \mB_{v}\boldsymbol{\Pi}^{\top},
\]
where $\bfR$ is a Jacobi smoother, $\boldsymbol{\Pi}$ is the N\'ed\'elec interpolation operator for vector functions, 
and $\mB_{\rm grad}$ and $\mB_{v}$ correspond to multigrid V-cycles for $\bfG^{\top}\bfS\bfG$ and $\boldsymbol{\Pi}^{\top}\bfS\boldsymbol{\Pi}$, respectively.  
In particular since $\bfC \bfG = \bf0$, we have 
\[
\bfG^{\top}\bfS\bfG 
= \bfG^{\top}\bfD^{\top}(\bfM_n^{\beta})^{-1}\bfD\bfG
= (\bfG^{\top}\bfM_e\bfG)^{\top} (\bfM_n^{\beta})^{-1} (\bfG^{\top}\bfM_e\bfG),
\]
which represents a fourth-order Laplacian operator, independent of $\bfsigma$ and unchanged through iterations,
and %$\boldsymbol{\Pi}^{\top}\bfS\boldsymbol{\Pi}: \bfH_0^1(\Omega) \to \bfH_0^1(\Omega)$ is the interpolation of $\bfS$ onto the smooth vector subspace. 
$\boldsymbol{\Pi}^{\top}\bfS\boldsymbol{\Pi}: \bfH_0^1(\Omega) \to \bfH_0^1(\Omega)$ is basically a Hodge Laplacian with a $\bfsigma$-dependent coefficient.
For $\bfG^{\top}\bfS\bfG$, we employ a standard geometric multigrid methods for Laplacians as its preconditioner,
while for $\boldsymbol{\Pi}^{\top}\bfS\boldsymbol{\Pi}$,
we exploit the multigrid methods for Hodge Laplacians.

\subsection{Numerical Results}

For the first example, we consider the nonlinear electromagnetic permeability function
\begin{equation*}
\nu(s) = a_0 + a_1 \exp(-a_2 s),
\end{equation*}
on the domain $\Omega = [-1,1]^3$. 
The exact solution is chosen as a trigonometric function:
\begin{equation*}
u(x_1, x_2, x_3) = \left [0, 0, \cos(\omega x_1) \cos(\omega x_2) \cos(\omega x_3)\right ]^{\top}
%\begin{bmatrix}
%0 \\[2pt]
%0 \\[2pt]
%1
%\end{bmatrix},
\end{equation*}
where $\omega$ is the frequency.

Note that $\Phi(s) = \nu(s)s$ and $\Phi^{\prime}(s) = a_0 + a_1(1 - a_2 s)\exp(-a_2 s)$. 
The monotonicity constant of $\Phi(\cdot)$ is $\underline{\nu} = a_0 - a_1 e^{-2}$, obtained by minimizing $\Phi^{\prime}(\cdot)$. 
We set $a_0 = 10$, $a_1 = 73.89$, and $a_2 = 1$. 
In this case, $\nu(s)s$ is strictly monotone but nearly not strongly monotone. 

We choose a fixed step size $\alpha = 0.5$ and the PCG tolerance ${\rm tol}_{\rm pcg} = 10^{-1}$ 
% \mnote{ is it too large? $10^{-2}$? Form HX preconditioner will take consierdable time. \JW{HX preconditioner is used as preconditioener for CG, so it is resembled every CG step. The time will note be reduced even epsilon decreases.}}
, meaning the relative residual is reduced below $10^{-1}$ at each PCG step. This corresponds to a quite rough approximation in the inner iteration. 
The iterative method terminates once the global relative residual decreases by a factor of $10^{-6}$.

The numerical errors with respect to mesh size $h$ and degrees of freedom (DoF) are reported in \cref{table: Maxwell numeircal error}. 
Here, $(u, \phi, \sigma)$ denote the exact solutions, and
$(u_h, \phi_h, \sigma_h)$ denote the numerical solutions obtained by the proposed algorithm. 
The observed convergence rates of the errors are at least first order. 
These results confirm that the proposed algorithm achieves good convergence behavior consistent with the expected accuracy of the finite element discretization. Next, we examine the performance of the solvers. 

% In \cref{fig:conver order}, we show that 
\begin{table}[ht]
  \centering
  \caption{Numerical error for nonlinear ferromagnetism model.}
  \label{table: Maxwell numeircal error}
  \renewcommand{\arraystretch}{1.25}
  \resizebox{.85\textwidth}{!}{
  \begin{tabular}{@{} c c c c c c c @{}}
    \toprule
$h$ & DoF	&	$\|u - u_h\|_{L^2}$	&	$\|\curl (u - u_h)\|_{L^2}$	&	$\|\sigma - \sigma_h\|_{L^2}$	&	$\| \phi - \phi_h\|_{L^2}$	\\ \hline
1/4 & 14129	&	1.17E+00	&	3.00E+00	&	3.64E+01	&	5.50E+00	\\
1/8 & 109665	&	4.28E-01	&	1.45E+00	&	1.76E+01	&	2.34E+00	\\
1/16 & 864449	&	1.39E-01	&	7.69E-01	&	9.75E+00	&	2.41E-01	\\
1/32 & 6865281	&	7.10E-02	&	3.98E-01	&	4.91E+00	&	1.50E-01	\\
1/64 & 54723329	&	3.50E-02	&	1.95E-01	&	2.44E+00	&	5.68E-02	\\
    \bottomrule
  \end{tabular}}
\end{table}

%\begin{figure}
%    \centering
%    \includegraphics[width=0.45\linewidth]{Figures/uorder.png}
%    \includegraphics[width=0.45\linewidth]{Figures/sigmaporder.png}
%    \caption{Order of convergence for the numerical error by the presented algorithm on nonlinear ferromagnetism model.}
%    \label{fig:conver order}
%\end{figure}

\begin{figure}[htbp]
    \centering
\begin{minipage}[t]{0.55\textwidth}
  \centering
  \captionof{table}{Computation performance for the proposed DualTPD algorithm on nonlinear ferromagnetism model.}
  \renewcommand{\arraystretch}{1.125}
%   \resizebox{.95\textwidth}{!}{\begin{tabular}{@{} c c c c c@{}}
%     \toprule
% $h$ 	&	 Iteration	&	\makecell[c]{PCG \\ (MG V-cycle)}	&	\makecell[c]{CPU time \\ (seconds)}	\\ \hline
% 1/4	&	18	&	141	&	2.6	\\
% 1/8	&	17	&	157	&	11	\\
% 1/16	&	16	&	149	&	103	\\
% 1/32	&	15	&	132	&	778	\\
% 1/64	&	14	&	138	&	6693	\\
%     \bottomrule
%   \end{tabular}}
    \resizebox{.95\textwidth}{!}{\begin{tabular}{@{} c c c c c@{}}
    \toprule
$h$ 	&	 Iteration	&	\makecell[c]{Avg. PCG \\ per step}	&	\makecell[c]{CPU time \\ in seconds}	\\ \hline
1/4	&	18	&	7.8&	2.6	\\
1/8	&	17	&	9.2	&	11	\\
1/16	&	16	&	9.3	&	103	\\
1/32	&	15	&	8.8	&	778	\\
1/64	&	14	&	9.9	&	6693	\\
    \bottomrule
  \end{tabular}}
  \label{table:dual-Maxwell-result}
 \end{minipage}
\hfill
 \begin{minipage}[t]{0.38\textwidth}
      \centering
          \captionof{figure}{Rate of time growth for DualTPD algorithm to solve nonlinear Maxwell equations.}
        \includegraphics[width=1\linewidth]{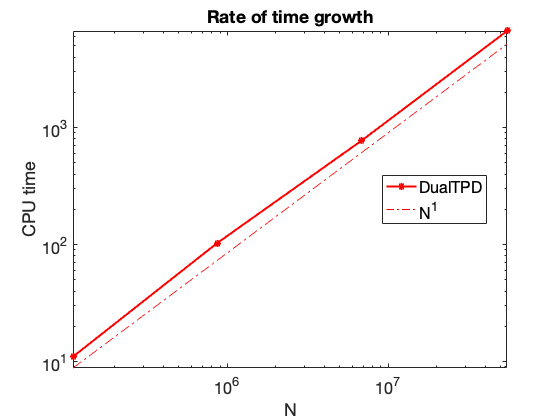}
    \label{fig:timegrowth_maxwell}
 \end{minipage}
 \end{figure}

We report the number of iterations and CPU time with respect to the mesh size 
$h$ in \cref{table:dual-Maxwell-result}. The most computationally expensive component is the PCG iteration. Therefore, we also include the averaged PCG step per iteration. 
% On average, approximately 8-9 PCG iterations are performed for one preconditioning step. 
In each PCG iteration, a single V-cycle multigrid is applied as a preconditioner. The iteration count of the proposed algorithm remains stable and independent of the mesh size $h$. In \cref{fig:timegrowth_maxwell}, we plot the rate of time growth with respect to the degree of freedom $(N)$. With effective preconditioning, the total computational time increases linearly with the problem size, demonstrating the good scalability of our algorithm.

% alpha = 0.5
%   \begin{minipage}[t]{0.95\linewidth}
%   \centering
%   \captionof{table}{Computation performance for the proposed algorithm. }
% \renewcommand{\arraystretch}{1.25}
% \begin{tabular}{@{} c c c c c@{}}
%     \toprule
% h & DoF	&	Iteration	&	PCG(Multigrid)	&	CPU Time	\\ \hline
% 1/4 & 13400	&	21	&	168	&	1.5s	\\
% 1/8 & 109665	&	21	&	165	&	6.5s	\\
% 1/16 & 828512	&	25	&	199	&	62s	\\
% 1/32 & 6590656	&	25	&	236	&	844s	\\
%     \bottomrule
% \end{tabular}
% \label{table: dual Newton result}
%   \end{minipage}

% \begin{figure}[htbp]
%     \centering
% 	\begin{minipage}{0.7\linewidth}
% 		\centering
% 		\vspace{-0.6cm}
% 		\setlength{\abovecaptionskip}{0.28cm}
% 		\includegraphics[width=\linewidth]{MPC8.pdf}
% 		\caption{The crane}
% 		\label{fig:1}
% 	\end{minipage}
% 	%\qquad
% 	\hfill
% 	\begin{minipage}{0.45\linewidth}
% 		\centering
%         \begin{tabular}{|c|c|}
%         \hline
%         aa & bb \\ \hline
%         cc & dd \\ \hline
%         \end{tabular}
%         \captionof{table}{The crane}%这里必须写table，不然标题就自动设置成figure
%         \label{fig:NormalSafe}
% 	    \end{minipage}
% \end{figure}

  \vspace{0.3cm}

We now compare our algorithm with other iterative methods. First, we consider a variant of the Chambolle--Pock algorithm~\cite{chambolle2011first} within the primal--dual framework (denoted DualPD). As we mentioned before, without efficient preconditioners, primal--dual methods are slow. So we also adapt our preconditioners below. For $\theta \in [0,1]$, given $(\bfsigma_k, \bfphi_k, \bar{\bfu}_k, \bfu_k)$, the updates are
\begin{equation}\label{eq:dualPD}
    \begin{aligned}
     \bfsigma_{k+1}  & = \bfsigma_k - \alpha_k \big( \bfJ(\bfM_{T}^{\gamma(\bfsigma_k)}\bfsigma_k) \big)^{-1} 
     \big(\bfM_{3T}^{\gamma(\bfsigma_k)} \bfsigma_k - \tilde \bfC \bar \bfu_k\big), \\
     \bfphi_{k+1}  & = \bfphi_k - \alpha_k (\bfM_n^{\beta})^{-1} \big(\bfM_n^{\beta}\bfphi_k - \bfD \bar \bfu_k\big), \\
     \bfu_{k+1} &=  \bfu_{k} - \alpha_k \, {\rm PCG} (\bfS_k; {\rm tol}_{\rm pcg}) \,  
     \big(\bfC^{\top} \bfsigma_{k+1} + \bfD^{\top}\bfphi_{k+1} - \bfJ \big), \\
     \bar{\bfu}_{k+1} &=  \bfu_{k+1} + \theta \big( \bfu_{k+1} -  \bfu_{k}\big).
    \end{aligned}
\end{equation}
Compared with DualTPD, no transformation is applied; instead, the updates for $\bfu_{k+1}$ use the new iterates $\bfsigma_{k+1}$ and $\bfphi_{k+1}$, together with an auxiliary variable $\bar{\bfu}_k$. We set $\alpha_k = 0.5$ and $\theta = 0.8$.
The associated block-diagonal preconditioner is
\begin{equation*}
\mathcal B_k^{-1} = 
\begin{pmatrix}
    ( \bfJ(\bfM_{T}^{\gamma(\bfsigma_k)}\bfsigma_k) )^{-1} & & \\ 
    & (\bfM_n^{\beta})^{-1} & \\
    & & {\rm PCG}(\bfS_k;{\rm tol}_{\rm pcg}) 
\end{pmatrix}.
\end{equation*}

\begin{table}[ht]
  \centering
  \caption{Computation performance for the three algorithms: DualTPD, DualPD and PrimalTPD, on the nonlinear ferromagnetism model.
  It shows that the three factors--duality, transformation and preconditioners--integratively accelerate the convergence.  }
  \renewcommand{\arraystretch}{1.125}
  \resizebox{.65\textwidth}{!}{  
  \begin{tabular}{@{} c |c| c c c @{}}
    \toprule
		&	$h$	&	DualTPD	&	DualPD	&	PrimalTPD \cite{chen2024transformed}	\\ \hline
\multirowcell{ 3}{Iteration}		&	1/4	&	18	&	26	&	32	\\
		&	1/8	&	17	&	26	&	48	\\
		&	1/16	&	16	&	33	&	72	\\ \hline
\multirowcell{ 3}{CPU time \\ in seconds}		&	1/4	&	2.6	&	3.7	&	2.3	\\
		&	1/8	&	11	&	20	&	44	\\ 
		&	1/16	&	103	&	228	&	765	\\ 
    \bottomrule
  \end{tabular}}
  \label{table:dual-Maxwell-comparison-result}
\end{table}

Additionally, we introduce the transformed primal--dual method (denoted PrimalTPD) proposed in~\cite{chen2024transformed} to the primal formulation that is considered in \cite{chen2024transformed}, which already outperforms other classical iterative methods such as fixed point and projected gradient descent methods based on the primal formulation as shown in \cite[Section 4]{chen2024transformed}.
In \cref{table:dual-Maxwell-comparison-result}, we list the iteration number and CPU time for the considered algorithms. 
The proposed algorithm is more efficient in both iteration number and computation time. Compared with other algorithms, the iteration number of DualTPD, is independent of $h$, which showcases the advantage of considering the dual formulation and using the transformed primal-dual method.

%In this case, due to the coupling between the differential operator and the nonlinear coefficient, the preconditioners only employ the weighted mass matrix with nonlinear coefficients. 

%In \cref{table:dual-Maxwell-comparison-result}, we list the iteration number and CPU time for the considered algorithms. Noted that PrimalTPD outperforms other classical iterative methods such as fixed point iteration and projected gradient descent as shown in \cite[Section 4]{chen2024transformed}. The proposed algorithms in this work are more efficient in both iteration number and computation time. Compared with other algorithms, the iteration number of DualTPD, is independent of $h$, which showcases the advantage of considering the dual formulation and using the transformed primal-dual method. 

% For this example, we also try the weighted mass matrices as preconditioner, but then only small step size is allowed for convergence and the algorithms fail to converge in reasonable time.

%  \end{minipage}

For the second example, we consider the $p$-curl problem arising in nonlinear electromagnetics, magnetohydrodynamics, and non-Newtonian fluids, where the nonlinearity of the curl operator follows a power-law relation~\cite{wan2020posteriori,hichmani2025mixed,miranda2010p}. 
We show that our algorithm can be naturally extended to $p$-curl problems by combining the discussion in Section~\ref{sec:dual p-Laplacian}. 
Fast solvers for $p$-curl problems are rarely studied, and this extension highlights the generality of our approach.  
The nonlinear coefficient is
\[
\nu(|\curl u|) = |\curl u|^{p-2},
\]
with the corresponding inverse relation
\[
\Phi^{-1}(|\bfsigma|) = |\bfsigma|^{p-1}, 
\qquad 
\gamma(\bfsigma) = |\bfsigma|^{p^*-2}.
\]
Thus, the nonlinear term of the $p$-curl problem in the dual formulation coincides with that of the $p$-Laplacian problem.  

Accordingly, we employ the Jacobian-based preconditioner \eqref{eq: Jacobian precon} for $\bfsigma$ with regularization parameters $\lambda = 10^{-4}$ and $\epsilon_0 = 10^{-16}$.  The remaining preconditioners are the same as in Section~\ref{sec: preconditioning}:
\[
\bfI_{\phi}^{-1} = ( \bfM^{\beta}_n)^{-1},
\qquad 
\bfI_{u_k}^{-1} = {\rm PCG} (\bfS_k, {\rm tol}_{\rm pcg}),
\]
where
\[
\bfS_k = \tilde \bfC^{\top}\bfI_{\sigma_k}^{-1} \tilde \bfC + \bfD^{\top}(\bfM_n^{\beta})^{-1} \bfD.
\]

\begin{table}[ht]
  \centering
  \caption{Comparison of DualTPD and DualPD algorithms on the $p$-curl problem.}
  \renewcommand{\arraystretch}{1.125}
  \resizebox{0.725\textwidth}{!}{
  \begin{tabular}{@{} c | c  c | c c | c  c | c c @{}}
    \toprule
	&	\multicolumn{4}{c|}{$p =1.5$}			&					\multicolumn{4}{c}{$p =3$}							\\ \cline{2-9}
	&	\multicolumn{2}{c|}{DualTPD}			&	\multicolumn{2}{c|}{DualPD} 			&	\multicolumn{2}{c|}{DualTPD}			&	\multicolumn{2}{c}{DualPD} 			\\ \hline
$\alpha$	&	\multicolumn{2}{c|}{$ 0.8$}			&	\multicolumn{2}{c|}{$0.7$}			&	\multicolumn{2}{c|}{$0.8$}			&	\multicolumn{2}{c}{$0.6$}			\\  \hline
$h$	&	 Iter. 	&	Time 	&	 Iter. 	&	Time 	&	 Iter. 	&	Time 	&	 Iter. 	&	Time 	\\\hline
1/4	&	11	&	1.3	&	74	&	4.2	&	11	&	1.3	&	68	&	3.3	\\
1/8	&	10	&	5.7	&	121	&	35	&	12	&	7.3	&	113	&	28	\\
1/16	&	10	&	57	&	200	&	490	&	12	&	72	&	185	&	417	\\
1/32	&	10	&	491	&	320	&	8981	&	12	&	647	&	283	&	8375	\\
    \bottomrule
  \end{tabular}}
  \label{table:dual-vs-dualPD-pcurl}
\end{table}

We test the $p$-curl problem for $p = 1.5$ and $p = 3$. 
The exact solution is chosen as a quadratic polynomial:
\[
u(x_1, x_2, x_3) =  \sum_{i=1}^3 (x_i + 1)(1-x_i)[0, 0, 1]^{\top}.
%\begin{bmatrix}
%0 \\[3pt]
%0 \\[3pt]
%\displaystyle \sum_{i=1}^3 (x_i + 1)(1-x_i)
%\end{bmatrix}.
\]
We set the step size $\alpha_k = 0.8$, and the PCG tolerance ${\rm tol}_{\rm pcg} = 10^{-1}$, meaning the relative residual is reduced below $10^{-1}$ at each PCG step. The iterative method terminates once the global relative residual decreases by a factor of $10^{-6}$.
The iteration counts, average PCG steps per iteration, and CPU times are reported in \cref{table:dual-vs-dualPD-pcurl}.  
For both $p < 2$ and $p > 2$, the iteration numbers remain stable, demonstrating the robust convergence of our algorithm.  
For comparison, \cref{table:dual-vs-dualPD-pcurl} also presents results for $p = 1.5$ and $p = 3$ using DualTPD and DualPD.  
%The step sizes are set to $\alpha_k = 0.8$ for $p = 1.5$ and $\alpha_k = 0.6$ for $p = 3$.  
DualTPD consistently outperforms DualPD, requiring far fewer iterations and substantially less CPU time.  

\section{Conclusion}\label{sec:conclusion}

We have proposed a DualTPD algorithm for nonlinear PDEs based on the Fenchel--Rockafellar duality and transformed primal-dual techniques. 
The dual formulation decouples nonlinear terms from differential operators and enables efficient local treatment of nonlinearities and effective preconditioner design. 
%Within this framework, triangular preconditioning yields globally convergent iterations and recovers Newton-type local acceleration through Jacobian-based preconditioners.
Applications to $p$-Laplacian and nonlinear Maxwell equations, including $p$-curl problems, have confirmed the efficiency and robustness of the approach. 
In both cases, the proposed methods achieves mesh-independent convergence, optimal CPU scaling, and often outperformed primal-based methods.

Future work includes refining preconditioners for strongly nonlinear regimes, incorporating adaptive strategies for step-size and preconditioner updates, and extending the framework to coupled multiphysics systems. 
Overall, the DualTPD algorithm provides a flexible and efficient tool for solving a broad class of nonlinear saddle-point problems.

% \begin{algorithm}
% \caption{Build tree}
% \label{alg:buildtree}
% \begin{algorithmic}
% \STATE{Define $P:=T:=\{ \{1\},\ldots,\{d\}$\}}
% \WHILE{$\#P > 1$}
% \STATE{Choose $C^\prime\in\mathcal{C}_p(P)$ with $C^\prime := \operatorname{argmin}_{C\in\mathcal{C}_p(P)} \varrho(C)$}
% \STATE{Find an optimal partition tree $T_{C^\prime}$ }
% \STATE{Update $P := (P{\setminus} C^\prime) \cup \{ \bigcup_{t\in C^\prime} t \}$}
% \STATE{Update $T := T \cup \{ \bigcup_{t\in\tau} t : \tau\in T_{C^\prime}{\setminus} \mathcal{L}(T_{C^\prime})\}$}
% \ENDWHILE
% \RETURN $T$
% \end{algorithmic}
% \end{algorithm}

% \appendix
% \section{An example appendix} 

%\section*{Acknowledgments}

\bibliographystyle{siamplain}
\bibliography{references}

\begin{thebibliography}{10}

\bibitem{aragon2023effective}
{\sc A.~Aragón, J.~F. Bonder, and D.~Rubio}, {\em Effective numerical
  computation of $p(x)$–{Laplace} equations in {2D}}, International Journal
  of Computer Mathematics, 100 (2023), pp.~2111--2123,
  \url{https://doi.org/10.1080/00207160.2023.2263103}.

\bibitem{arnol2013mathematical}
{\sc V.~I. Arnol'd}, {\em Mathematical methods of classical mechanics},
  vol.~60, Springer Science \& Business Media, 2013,
  \url{https://doi.org/10.1007/978-1-4757-2063-1}.

\bibitem{1958ArrowHurwiczUzawa}
{\sc K.~Arrow, L.~Hurwicz, and H.~Uzawa}, {\em Studies in Linear and Nonlinear
  Programming}, Stanford University Press, 1958.

\bibitem{attouch2016strongly}
{\sc H.~Attouch, L.~M. Brice{\~n}o-Arias, and P.~L. Combettes}, {\em A strongly
  convergent primal--dual method for nonoverlapping domain decomposition},
  Numerische Mathematik, 133 (2016), pp.~443--470,
  \url{https://doi.org/10.1007/s00211-015-0751-4}.

\bibitem{1995DominiqueJeanPaul}
{\sc D.~Az\'e and J.-P. Penot}, {\em Uniformly convex and uniformly smooth
  convex functions}, Annales de la Facult\'e des sciences de Toulouse :
  Math\'ematiques, Ser. 6, 4 (1995), pp.~705--730,
  \url{https://www.numdam.org/item/AFST_1995_6_4_4_705_0/}.

\bibitem{bachinger2005numerical}
{\sc F.~Bachinger, U.~Langer, and J.~Sch{\"o}berl}, {\em Numerical analysis of
  nonlinear multiharmonic eddy current problems}, Numerische Mathematik, 100
  (2005), pp.~593--616, \url{https://doi.org/10.1007/s00211-005-0597-2}.

\bibitem{balseiro2020dual}
{\sc S.~Balseiro, H.~Lu, and V.~Mirrokni}, {\em Dual mirror descent for online
  allocation problems}, in International Conference on Machine Learning, PMLR,
  2020, pp.~613--628.

\bibitem{barrett1993finite}
{\sc J.~W. Barrett and W.~B. Liu}, {\em Finite element approximation of the
  $p$-{Laplacian}}, Mathematics of Computation, 61 (1993), pp.~523--537,
  \url{http://www.jstor.org/stable/2153239} (accessed 2025-09-09).

\bibitem{bartels2021error}
{\sc S.~Bartels}, {\em Error estimates for a class of discontinuous {Galerkin}
  methods for nonsmooth problems via convex duality relations}, Mathematics of
  Computation, 90 (2021), pp.~2579--2602,
  \url{https://doi.org/10.1090/mcom/3821}.

\bibitem{bartels2021nonconforming}
{\sc S.~Bartels}, {\em Nonconforming discretizations of convex minimization
  problems and precise relations to mixed methods}, Computers \& Mathematics
  with Applications, 93 (2021), pp.~214--229.

\bibitem{benzi2005numerical}
{\sc M.~Benzi, G.~H. Golub, and J.~Liesen}, {\em Numerical solution of saddle
  point problems}, Acta Numerica, 14 (2005), p.~1–137,
  \url{https://doi.org/10.1017/S0962492904000212}.

\bibitem{benzi2011modified}
{\sc M.~Benzi, M.~A. Olshanskii, and Z.~Wang}, {\em Modified augmented
  {Lagrangian} preconditioners for the incompressible {Navier--Stokes}
  equations}, International Journal for Numerical Methods in Fluids, 66 (2011),
  pp.~486--508, \url{https://doi.org/10.1002/fld.2267}.

\bibitem{borwein2006convex}
{\sc J.~Borwein and A.~Lewis}, {\em Convex Analysis and Nonlinear Optimization:
  Theory and Examples}, CMS Books in Mathematics, Springer New York, 2005.

\bibitem{bramble1988preconditioning}
{\sc J.~H. Bramble and J.~E. Pasciak}, {\em A preconditioning technique for
  indefinite systems resulting from mixed approximations of elliptic problems},
  Mathematics of Computation, 50 (1988), pp.~1--17.

\bibitem{brandt1977multi}
{\sc A.~Brandt}, {\em Multi-level adaptive solutions to boundary-value
  problems}, Mathematics of Computation, 31 (1977), pp.~333--390,
  \url{https://doi.org/10.1090/S0025-5718-1977-0431719-X}.

\bibitem{breit2015finite}
{\sc D.~Breit, L.~Diening, and S.~Schwarzacher}, {\em Finite element
  approximation of the $p(\cdot)$-{Laplacian}}, SIAM Journal on Numerical
  Analysis, 53 (2015), pp.~551--572, \url{https://doi.org/10.1137/130946046}.

\bibitem{briceno2011monotone+}
{\sc L.~M. Brice\~{n}o Arias and P.~L. Combettes}, {\em A monotone+skew
  splitting model for composite monotone inclusions in duality}, SIAM Journal
  on Optimization, 21 (2011), pp.~1230--1250,
  \url{https://doi.org/10.1137/10081602X}.

\bibitem{brune2015composing}
{\sc P.~R. Brune, M.~G. Knepley, B.~F. Smith, and X.~Tu}, {\em Composing
  scalable nonlinear algebraic solvers}, SIAM Review, 57 (2015), pp.~535--565,
  \url{https://doi.org/10.1137/130936725}.

\bibitem{carrillo2022primal}
{\sc J.~A. Carrillo, K.~Craig, L.~Wang, and C.~Wei}, {\em Primal dual methods
  for {Wasserstein} gradient flows}, Foundations of Computational Mathematics,
  (2022), pp.~1--55.

\bibitem{carstensen2012mixed}
{\sc C.~Carstensen, D.~G\"{u}nther, and H.~Rabus}, {\em Mixed finite element
  method for a degenerate convex variational problem from topology
  optimization}, SIAM Journal on Numerical Analysis, 50 (2012), pp.~522--543,
  \url{https://doi.org/10.1137/100806837}.

\bibitem{carstensen2015nonconforming}
{\sc C.~Carstensen and D.~J. Liu}, {\em Nonconforming {FEMs} for an optimal
  design problem}, SIAM Journal on Numerical Analysis, 53 (2015), pp.~874--894,
  \url{https://doi.org/10.1137/130927103}.

\bibitem{carstensen2021unstabilized}
{\sc C.~Carstensen and T.~Tran}, {\em Unstabilized hybrid high-order method for
  a class of degenerate convex minimization problems}, SIAM Journal on
  Numerical Analysis, 59 (2021), pp.~1348--1373,
  \url{https://doi.org/10.1137/20M1335625}.

\bibitem{chambolle2011first}
{\sc A.~Chambolle and T.~Pock}, {\em A first-order primal-dual algorithm for
  convex problems with applications to imaging}, Journal of Mathematical
  Imaging and Vision, 40 (2011), pp.~120--145.

\bibitem{chen2024transformed}
{\sc L.~Chen, R.~Guo, and J.~Wei}, {\em Transformed primal-dual methods with
  variable-preconditioners}, 2023, \url{https://arxiv.org/abs/2312.12355}.

\bibitem{chen2023transformed}
{\sc L.~Chen and J.~Wei}, {\em Transformed primal--dual methods for nonlinear
  saddle point systems}, Journal of Numerical Mathematics, 31 (2023),
  pp.~281--311.

\bibitem{chen2018convergence}
{\sc L.~Chen and Y.~Wu}, {\em Convergence analysis for a class of iterative
  methods for solving saddle point systems}, 2018,
  \url{https://arxiv.org/abs/1710.03409},
  \url{https://arxiv.org/abs/1710.03409}.

\bibitem{chen1998global}
{\sc X.~Chen}, {\em Global and superlinear convergence of inexact {Uzawa}
  methods for saddle point problems with nondifferentiable mappings}, SIAM
  Journal on Numerical Analysis, 35 (1998), pp.~1130--1148,
  \url{https://doi.org/10.1137/S0036142995295789}.

\bibitem{chen1998preconditioned}
{\sc X.~Chen}, {\em On preconditioned {Uzawa} methods and {SOR} methods for
  saddle-point problems}, Journal of Computational and Applied Mathematics, 100
  (1998), pp.~207--224, \url{https://doi.org/10.1016/S0377-0427(98)00197-6}.

\bibitem{clason2017primal}
{\sc C.~Clason and T.~Valkonen}, {\em Primal-dual extragradient methods for
  nonlinear nonsmooth pde-constrained optimization}, SIAM Journal on
  Optimization, 27 (2017), pp.~1314--1339,
  \url{https://doi.org/10.1137/16M1080859}.

\bibitem{2007CombettesPesquet}
{\sc P.~L. Combettes and J.-C. Pesquet}, {\em A {D}ouglas--{R}achford splitting
  approach to nonsmooth convex variational signal recovery}, IEEE Journal of
  Selected Topics in Signal Processing, 1 (2007), pp.~564--574,
  \url{https://doi.org/10.1109/JSTSP.2007.910264}.

\bibitem{combettes2005signal}
{\sc P.~L. Combettes and V.~R. Wajs}, {\em Signal recovery by proximal
  forward-backward splitting}, Multiscale Modeling \& Simulation, 4 (2005),
  pp.~1168--1200, \url{https://doi.org/10.1137/050626090}.

\bibitem{creuse2007posteriori}
{\sc E.~Creuse, M.~Farhloul, and L.~Paquet}, {\em A posteriori error estimation
  for the dual mixed finite element method for the $p$-{Laplacian} in a
  polygonal domain}, Computer Methods in Applied Mechanics and Engineering, 196
  (2007), pp.~2570--2582, \url{https://doi.org/10.1016/j.cma.2006.11.023}.

\bibitem{ekeland1999convex}
{\sc I.~Ekeland and R.~Témam}, {\em Convex Analysis and Variational Problems},
  Society for Industrial and Applied Mathematics, 1999,
  \url{https://doi.org/10.1137/1.9781611971088}.

\bibitem{1997Evans}
{\sc L.~Evans}, {\em Partial Differential Equations}, Graduate Studies in
  Mathematics, American Mathematical Society, 2022.

\bibitem{fang2014single}
{\sc F.~Fang, F.~Li, and T.~Zeng}, {\em Single image dehazing and denoising:
  {A} fast variational approach}, SIAM Journal on Imaging Sciences, 7 (2014),
  pp.~969--996, \url{https://doi.org/10.1137/130919696}.

\bibitem{farhloul2000mixed}
{\sc M.~Farhloul and H.~Manouzi}, {\em On a mixed finite element method for the
  $p$-{Laplacian}}, Canadian Applied Mathematics Quarterly, 8 (2000),
  pp.~67--78.

\bibitem{fortin2000augmented}
{\sc M.~Fortin and R.~Glowinski}, {\em Augmented Lagrangian Methods:
  Applications to the Numerical Solution of Boundary-Value Problems}, Studies
  in Mathematics and its Applications, North Holland, 2000.

\bibitem{gabay1983chapter}
{\sc D.~Gabay}, {\em Chapter {IX} {Applications} of the method of multipliers
  to variational inequalities}, in Augmented Lagrangian Methods: Applications
  to the Numerical Solution of Boundary-Value Problems, M.~Fortin and
  R.~Glowinski, eds., vol.~15 of Studies in Mathematics and Its Applications,
  Elsevier, 1983, pp.~299--331,
  \url{https://doi.org/10.1016/S0168-2024(08)70034-1}.

\bibitem{gabay1976dual}
{\sc D.~Gabay and B.~Mercier}, {\em A dual algorithm for the solution of
  nonlinear variational problems via finite element approximation}, Computers
  \& Mathematics with Applications, 2 (1976), pp.~17--40,
  \url{https://doi.org/10.1016/0898-1221(76)90003-1}.

\bibitem{goldstein2009split}
{\sc T.~Goldstein and S.~Osher}, {\em The split {Bregman} method for
  {L}1-regularized problems}, SIAM Journal on Imaging Sciences, 2 (2009),
  pp.~323--343, \url{https://doi.org/10.1137/080725891}.

\bibitem{hestenes1969multiplier}
{\sc M.~R. Hestenes}, {\em Multiplier and gradient methods}, Journal of
  Optimization Theory and Applications, 4 (1969), pp.~303--320,
  \url{https://doi.org/10.1007/BF00927673}.

\bibitem{hichmani2025mixed}
{\sc M.~Hichmani, M.~Laforest, and E.~M. Zaoui}, {\em Mixed formulation in {3D}
  for a nonlinear eddy current problem from applied superconductivity},
  Discrete and Continuous Dynamical Systems, 46 (2026), pp.~287--304,
  \url{https://doi.org/10.3934/dcds.2025100}.

\bibitem{hiptmair2007nodal}
{\sc R.~Hiptmair and J.~Xu}, {\em Nodal auxiliary space preconditioning in
  {H}(curl) and {H}(div) spaces}, SIAM Journal on Numerical Analysis, 45
  (2007), pp.~2483--2509, \url{https://doi.org/10.1137/060660588}.

\bibitem{hu2001iterative}
{\sc Q.~Hu and J.~Zou}, {\em An iterative method with variable relaxation
  parameters for saddle-point problems}, SIAM Journal on Matrix Analysis and
  Applications, 23 (2001), pp.~317--338,
  \url{https://doi.org/10.1137/S0895479899364064}.

\bibitem{hu2002two}
{\sc Q.~Hu and J.~Zou}, {\em Two new variants of nonlinear inexact {Uzawa}
  algorithms for saddle-point problems}, Numerische Mathematik, 93 (2002),
  pp.~333--359, \url{https://doi.org/10.1007/s002110100386}.

\bibitem{hu2006nonlinear}
{\sc Q.~Hu and J.~Zou}, {\em Nonlinear inexact {Uzawa} algorithms for linear
  and nonlinear saddle-point problems}, SIAM Journal on Optimization, 16
  (2006), pp.~798--825, \url{https://doi.org/10.1137/S1052623403428683}.

\bibitem{huang2007preconditioned}
{\sc Y.~Q. Huang, R.~Li, and W.~Liu}, {\em Preconditioned descent algorithms
  for $p$-{Laplacian}}, Journal of Scientific Computing, 32 (2007),
  pp.~343--371, \url{https://doi.org/10.1007/s10915-007-9134-z}.

\bibitem{kolev2009parallel}
{\sc T.~V. Kolev and P.~S. Vassilevski}, {\em Parallel auxiliary space {AMG}
  for {H}(curl) problems}, Journal of Computational Mathematics, 27 (2009),
  pp.~604--623.

\bibitem{landau1960mechanics}
{\sc L.~Landau and E.~Lifshitz}, {\em Mechanics: Volume 1}, Course of
  Theoretical Physics, Butterworth-Heinemann, 1976.

\bibitem{Lindqvist2019}
{\sc P.~Lindqvist}, {\em Notes on the Stationary $p$-{Laplace} Equation},
  SpringerBriefs in Mathematics, Springer Cham, 2019,
  \url{https://doi.org/10.1007/978-3-030-14501-9}.

\bibitem{liu2023fast}
{\sc H.~Liu and D.~Wang}, {\em Fast operator splitting methods for obstacle
  problems}, Journal of Computational Physics, 477 (2023), p.~111941,
  \url{https://doi.org/10.1016/j.jcp.2023.111941}.

\bibitem{miranda2010p}
{\sc F.~Miranda, J.-F. Rodrigues, and L.~Santos}, {\em On a $p$-curl system
  arising in electromagnetism}, 2010, \url{https://arxiv.org/abs/1009.0424}.

\bibitem{neuberger2013newton}
{\sc J.~M. Neuberger, N.~Sieben, and J.~W. Swift}, {\em {Newton}'s method and
  symmetry for semilinear elliptic {PDE} on the cube}, SIAM Journal on Applied
  Dynamical Systems, 12 (2013), pp.~1237--1279,
  \url{https://doi.org/10.1137/120899054}.

\bibitem{pollock2015regularized}
{\sc S.~Pollock}, {\em A regularized {Newton}-like method for nonlinear {PDE}},
  Numerical Functional Analysis and Optimization, 36 (2015), pp.~1493--1511,
  \url{https://doi.org/10.1080/01630563.2015.1069328}.

\bibitem{potter1993dual}
{\sc L.~C. Potter and K.~S. Arun}, {\em A dual approach to linear inverse
  problems with convex constraints}, SIAM Journal on Control and Optimization,
  31 (1993), pp.~1080--1092, \url{https://doi.org/10.1137/0331049}.

\bibitem{rockafellar1997convex}
{\sc R.~Rockafellar}, {\em Convex Analysis}, Princeton Landmarks in Mathematics
  and Physics, Princeton University Press, 1997.

\bibitem{rockafellar1976augmented}
{\sc R.~T. Rockafellar}, {\em Augmented {Lagrangians} and applications of the
  proximal point algorithm in convex programming}, Mathematics of Operations
  Research, 1 (1976), pp.~97--116, \url{https://doi.org/10.1287/moor.1.2.97}.

\bibitem{1977Scheurer}
{\sc B.~Scheurer}, {\em Existence et approximation de points selles pour
  certains probl\`emes non lin\'eaires}, RAIRO. Analyse num\'erique, 11 (1977),
  pp.~369--400, \url{https://www.numdam.org/item/M2AN_1977__11_4_369_0/}.

\bibitem{sidky2012convex}
{\sc E.~Y. Sidky, J.~H. Jørgensen, and X.~Pan}, {\em Convex optimization
  problem prototyping for image reconstruction in computed tomography with the
  {Chambolle}–{Pock} algorithm}, Physics in Medicine \& Biology, 57 (2012),
  p.~3065, \url{https://doi.org/10.1088/0031-9155/57/10/3065}.

\bibitem{song2019inexact}
{\sc Y.~Song, X.~Yuan, and H.~Yue}, {\em An inexact {Uzawa} algorithmic
  framework for nonlinear saddle point problems with applications to elliptic
  optimal control problem}, SIAM Journal on Numerical Analysis, 57 (2019),
  pp.~2656--2684, \url{https://doi.org/10.1137/19M1245736}.

\bibitem{tran2024discrete}
{\sc N.~T. Tran}, {\em Discrete weak duality of hybrid high-order methods for
  convex minimization problems}, SIAM Journal on Numerical Analysis, 62 (2024),
  pp.~1492--1514, \url{https://doi.org/10.1137/23M1594534}.

\bibitem{valkonen2014primal}
{\sc T.~Valkonen}, {\em A primal–dual hybrid gradient method for nonlinear
  operators with applications to {MRI}}, Inverse Problems, 30 (2014),
  p.~055012, \url{https://doi.org/10.1088/0266-5611/30/5/055012},
  \url{https://dx.doi.org/10.1088/0266-5611/30/5/055012}.

\bibitem{wan2020posteriori}
{\sc A.~T.~S. Wan and M.~Laforest}, {\em A posteriori error estimation for the
  $p$-curl problem}, SIAM Journal on Numerical Analysis, 58 (2020),
  pp.~460--491, \url{https://doi.org/10.1137/16M1075624}.

\bibitem{xu2020adaptive}
{\sc Y.~Xu, I.~Yousept, and J.~Zou}, {\em An adaptive edge element
  approximation of a quasilinear {H}(curl)-elliptic problem}, Mathematical
  Models and Methods in Applied Sciences, 30 (2020), pp.~2799--2826,
  \url{https://doi.org/10.1142/S0218202520500554}.

\bibitem{yousept2013optimal}
{\sc I.~Yousept}, {\em Optimal control of quasilinear
  $\boldsymbol{H}(\mathbf{curl})$-elliptic partial differential equations in
  magnetostatic field problems}, SIAM Journal on Control and Optimization, 51
  (2013), pp.~3624--3651, \url{https://doi.org/10.1137/120904299}.

\bibitem{yu2011dual}
{\sc H.-F. Yu, F.-L. Huang, and C.-J. Lin}, {\em Dual coordinate descent
  methods for logistic regression and maximum entropy models}, Machine
  Learning, 85 (2011), pp.~41--75,
  \url{https://doi.org/10.1007/s10994-010-5221-8}.

\bibitem{zosso2017efficient}
{\sc D.~Zosso, B.~Osting, M.~Xia, and S.~J. Osher}, {\em An efficient
  primal-dual method for the obstacle problem}, Journal of Scientific
  Computing, 73 (2017), pp.~416--437,
  \url{https://doi.org/10.1007/s10915-017-0420-0}.

\bibitem{zulehner2002analysis}
{\sc W.~Zulehner}, {\em Analysis of iterative methods for saddle point
  problems: {A} unified approach}, Mathematics of Computation, 71 (2002),
  pp.~479--505, \url{https://doi.org/10.1090/S0025-5718-01-01324-2}.

\end{thebibliography}
\end{document}